\documentclass[review]{elsarticle}
\pdfoutput=1

\usepackage[table,xcdraw]{xcolor}  % xcolor for table coloring
\usepackage{tikz}  % TikZ for graphics
\usepackage{lineno,hyperref}  % Line numbering, hyperlinks
\usepackage{mathtools, amssymb, amsmath}  % Math packages
\usepackage{circuitikz, graphicx, epstopdf}  % Electrical circuits, graphics
\usetikzlibrary{circuits.ee.IEC, fadings, patterns, shadows.blur, shapes}  % Consolidate TikZ libraries

\usepackage{float}
\usepackage{natbib}  % Bibliography management
\usepackage{makecell}  % Advanced cell formatting in tables
\usepackage{multirow, tabularx, booktabs}  % Table formatting
\usepackage{extarrows}  % Extended arrows in equations
\usepackage{algorithm, algpseudocode}  % Algorithms
\usepackage{siunitx}  % Handling units
\usepackage{gensymb}  % Generic symbols like degrees

\usepackage{colortbl}  % Additional table color features
\usepackage{subcaption}

\journal{Computers in Biology and Medicine}
\usepackage{geometry}  % Custom page margins
\newgeometry{vmargin={15mm}, hmargin={12mm,17mm}}
\begin{document}

\begin{frontmatter}

\title{A comparative analysis of metamodels for 0D cardiovascular models, and pipeline for sensitivity analysis, parameter estimation, and uncertainty quantification}

\author[1]{John M. Hanna}
\author[1]{Pavlos Varsos}
\author[1]{Jérôme Kowalski}
\author[1,2]{Lorenzo Sala}
\author[1]{Roel Meiburg}
\author[1]{Irene E. Vignon-Clementel}

\address[1]{Inria, Research Center Saclay Ile-de-France, France}
\address[2]{Université Paris-Saclay, INRAE, MaIAGE, 78350, Jouy-en-Josas, France}

\begin{abstract}

Zero-dimensional (0D) cardiovascular models are reduced-order models aimed at studying the global dynamics of the whole circulation system or transport within it. They are employed to obtain estimates of important biomarkers for surgery planning and assessment applications (such as pressures, volumes, flow rates, and concentrations in the circulation) and can provide boundary conditions for high-fidelity three-dimensional models. Despite their low computational cost, tasks such as parameter estimation or uncertainty quantification require a large number of model evaluations, which is still a computationally expensive task. This motivates the building of metamodels in an offline stage, which can be evaluated significantly faster than 0D models. In this work, a pipeline going from 0D cardiovascular models to the building of metamodels and showcasing their use for tasks such as sensitivity analysis, parameter estimation, or uncertainty quantification is proposed. Three different strategies are assessed to build metamodels for 0D cardiovascular models, namely Neural Networks, Polynomial Chaos Expansion, and Gaussian Processes. The metamodels are assessed for three different 0D models. The first is a lumped model aimed at predicting the pressure in the portal vein after surgery. Due to the strong interaction between local liver hemodynamics and global circulation, the full circulation is modeled. The second one is simulating the whole-body circulation under the conditions of pulmonary arterial hypertension before and after shunt insertion. The final model is aimed at assessing the blood perfusion of an organ after a revascularization surgery. The transport of a contrast agent is modeled on top of a simplified 0D hemodynamics model. This model is chosen due to the different nature of the output which is a signal (concentration of the contrast agent over time), which requires a different treatment from the metamodeling point of view. The metamodels are trained and tested on synthetic data generated from the 0D models. It was found that neural networks offer the most convenient way of building metamodels in terms of the quality of the results, computational time, and practical ease of performing parameter estimation, sensitivity analysis, or uncertainty quantification tasks. Finally, we demonstrate a full pipeline of sensitivity analysis, inverse problem and (patient-specific) UQ, with a neural network as emulator.

\end{abstract}


\end{frontmatter}

%\linenumbers

\section{Introduction}
Many surgeries and percutaneous procedures, e.g., to treat cardiovascular diseases, organ cancer, or chronic diseases, lead to post-procedural changes in the blood circulation. These changes affect pressure and perfusion of organs and can lead to significant morbidity or mortality, but are often difficult to predict or assess based on pre-procedural metrics. Advances in the field of computational mechanics and increased computational power have made it possible to assist in personalized treatment planning or assessment by predicting previously mentioned hemodynamic changes after intervention. 

Three-dimensional (3D) models of blood flow dynamics and heart or vessel wall mechanics can provide high-fidelity solutions, for instance, using finite element \cite{fem} or finite volume methods \cite{fvm}. Additionally, less widely used methods, like meshless approaches, including the smoothed-particle hydrodynamics \cite{sph}, can also be adopted. In general, 3D simulations can be seen as a powerful tool for predicting clinically relevant biomarkers such as pressures, blood flow distribution, and wall shear stresses. However, the accuracy of 3D methods comes at the expense of high computational cost, which can reach days even on high-performance computers. This can be an issue when real-time predictions are needed, for instance, during surgery. It is also an issue for surgery optimization and sensitivity analysis problems when several simulations are required. Furthermore, these 3D models are often only feasible when modeling a specific vessel and/or organ, leaving full 3D simulations of the entire body circulation out of reach.

Reduced order models offer an attractive approach to alleviate this issue. Zero-dimensional (0D) cardiovascular models are developed to simulate the whole circulation hemodynamics \cite{0D-1D_review, caiazzo:hal-01955520}. 0D models are based on the analogy between hydraulic and electrical systems, connecting different components of the circulation. This analogy leads to a system of differential-algebraic equations (DAE, coupled ordinary differential equations and algebraic equations), representing mass conservation, momentum balance and constitutive laws in different circulatory compartments. 0D models provide faster predictions of pressures, volumes, flow rates, and tracer concentration in the whole cardiovascular system,  compared to 3D models since they only depend on solving a set of DAEs. This leads to reduced accuracy due to the simplification of the highly complex system and geometry. They can be stand-alone or provide boundary conditions for high fidelity 3D models \cite{pant2022multiscale}.

However, for clinical applications, the system needs to be solved hundreds to hundreds of thousands of times to solve complex \emph{meta-problems}, including sensitivity analysis (SA), parameter estimation, or uncertainty quantification (UQ). This can usually take hours or days in certain cases. To address this, data-driven \emph{meta-models} -- or emulators -- can be built offline from a dataset generated by running the 0D model multiple times and varying the input parameters. After such a model is built, it can be run online for even faster predictions than the 0D model itself. Moreover, the surrogate model is usually differentiable, which can then be utilized to solve SA problems or parameter estimation in a faster manner.

Various methods can be employed to build metamodels, with feed-forward neural networks (FFNN) being one of the most widely used approaches across different applications. FFNNs are favored for their ability to effectively approximate any continuous function when the appropriate weights (network parameters) are applied, a concept supported by the universal approximation theorem \cite{universal_approx}. Other neural network architectures have also demonstrated significant success in specific domains; for example, convolutional neural networks (CNNs) excel in computer vision tasks \cite{computer_vision}, while recurrent neural networks (RNNs) are highly effective for time series data, particularly in natural language processing \cite{nlp}. Neural networks generally perform exceptionally well when large datasets are available. However, in the context of building metamodels for hemodynamics, comprehensive and extensive patient data are often lacking. As an alternative, this database can be generated from 0D models. Despite the widespread application of neural networks in other fields, their use in developing reduced cardiovascular models is still in the exploratory stages.

Despite the increased interest in the topic, few research works have been published on the application of neural networks for 0D cardiovascular models. Bonnemain et al. predicted the left ventricular (LV) systolic function from systemic and pulmonary arterial pressure signals by neural networks \cite{LV_nn, LV_device_2021}. Significant efforts have also been made to create surrogates of three-dimensional models to reduce computational costs by exploiting various techniques; these include Graph Neural Networks and Physics-Informed Neural Networks (PINNs) \cite{kissas2020machine, de2024finite, pegolotti2024learning} or Latent Neural Networks \cite{salvador2024whole}.

Other methods have also been employed over the years, such as polynomial chaos expansion (PCE). Originally developed to represent a random output variable in terms of a polynomial function of multiple random input variables, \cite{wiener1938homogeneous} PCE exploits the properties of orthogonal polynomials to generate a surrogate model. It is mainly popular in the context of SA, as it can be used to approximate the decomposition of variance of model output with respect to model input, thereby generating Sobol sensitivity indices \cite{sudret2008global, sala2023sensitivity}. While PCE has been used in cardiovascular applications, it has seen more usage  in sensitivity analysis \cite{eck2016guide, meiburg2020uncertainty, schaferglobal}  than as an emulator \cite{heinen2020metamodeling}.

Another type of methodology is the Gaussian Processes (GPs). 
GPs are a versatile non-parametric statistical framework used for regression and classification tasks, providing a systematic way to quantify uncertainty. This makes them highly useful across a wide range of fields, for instance, in the field of machine learning \cite{seeger2004gaussian}, robotics \cite{deisenroth2013gaussian}, bioinformatics \cite{lawrence2005probabilistic}, and medicine \cite{coveney2020sensitivity}. In the field of cardiovascular models, GPs have been applied to the study of ventricular mechanics \cite{di2018gaussian, strocchi2023cell} and pulmonary circulation \cite{paun2021markov}, while more recently, they have been implemented for reconstructing blood flow regimes in realistic geometries \cite{ashtiani2024reconstructing}.

In this paper, a pipeline to build surrogate models for 0D cardiovascular systems is proposed. FFNN, PCE, and GPs are assessed to build surrogates for three different cardiovascular models. The first model is a simplified whole body hemodynamics model which predicts portal vein pressure following liver surgery \cite{audebert2017partial,golse2021predicting}. The second one is a more complex whole body model which aids in intervention planning for Pulmonary Arterial Hypertension (PAH) patients \cite{pant2022multiscale}. Finally, the third model evaluates blood perfusion in an organ post-revascularization surgery via its effect on the concentration in time of an injected tracer. The first two models were selected based on their complexity, as they simulate the whole-body circulation under different clinical scenarios and generate scalar biomarkers as outputs. In contrast, the third model, though simplified, was chosen because its output of interest is a time-dependent signal, requiring a different method to build the surrogate model. More details about the three models are discussed in section~\ref{section:0D}. The application of surrogate models for conducting SA, solving inverse problems and performing UQ is also investigated as these steps are crucial to create the digital twin of a patient.

\begin{figure}[H]
    \centering
    \resizebox{0.8\textwidth}{!}{\tikzset{every picture/.style={line width=0.75pt}} %set default line width to 0.75pt        

\begin{tikzpicture}[x=0.75pt,y=0.75pt,yscale=-1,xscale=1]
%uncomment if require: \path (0,398); %set diagram left start at 0, and has height of 398

%Rounded Rect [id:dp7195095635185288] 
\draw  [fill={rgb, 255:red, 155; green, 155; blue, 155 }  ,fill opacity=1 ] (270,18) .. controls (270,13.58) and (273.58,10) .. (278,10) -- (402,10) .. controls (406.42,10) and (410,13.58) .. (410,18) -- (410,42) .. controls (410,46.42) and (406.42,50) .. (402,50) -- (278,50) .. controls (273.58,50) and (270,46.42) .. (270,42) -- cycle ;

%Rounded Rect [id:dp8473981060647614] 
\draw  [fill={rgb, 255:red, 155; green, 155; blue, 155 }  ,fill opacity=1 ] (270,105.5) .. controls (270,96.94) and (276.94,90) .. (285.5,90) -- (394.5,90) .. controls (403.06,90) and (410,96.94) .. (410,105.5) -- (410,134.5) .. controls (410,143.06) and (403.06,150) .. (394.5,150) -- (285.5,150) .. controls (276.94,150) and (270,143.06) .. (270,134.5) -- cycle ;
%Rounded Rect [id:dp0659761278353116] 
\draw  [fill={rgb, 255:red, 80; green, 227; blue, 194 }  ,fill opacity=1 ] (10,18) .. controls (10,13.58) and (13.58,10) .. (18,10) -- (142,10) .. controls (146.42,10) and (150,13.58) .. (150,18) -- (150,42) .. controls (150,46.42) and (146.42,50) .. (142,50) -- (18,50) .. controls (13.58,50) and (10,46.42) .. (10,42) -- cycle ;

%Rounded Rect [id:dp12218856881078399] 
\draw  [fill={rgb, 255:red, 245; green, 166; blue, 35 }  ,fill opacity=1 ] (540,18) .. controls (540,13.58) and (543.58,10) .. (548,10) -- (672,10) .. controls (676.42,10) and (680,13.58) .. (680,18) -- (680,42) .. controls (680,46.42) and (676.42,50) .. (672,50) -- (548,50) .. controls (543.58,50) and (540,46.42) .. (540,42) -- cycle ;

%Rounded Rect [id:dp6946918690723931] 
\draw  [fill={rgb, 255:red, 80; green, 227; blue, 194 }  ,fill opacity=1 ] (10,108) .. controls (10,103.58) and (13.58,100) .. (18,100) -- (142,100) .. controls (146.42,100) and (150,103.58) .. (150,108) -- (150,132) .. controls (150,136.42) and (146.42,140) .. (142,140) -- (18,140) .. controls (13.58,140) and (10,136.42) .. (10,132) -- cycle ;

%Rounded Rect [id:dp10569131010431665] 
\draw  [fill={rgb, 255:red, 245; green, 166; blue, 35 }  ,fill opacity=1 ] (540,108) .. controls (540,103.58) and (543.58,100) .. (548,100) -- (672,100) .. controls (676.42,100) and (680,103.58) .. (680,108) -- (680,132) .. controls (680,136.42) and (676.42,140) .. (672,140) -- (548,140) .. controls (543.58,140) and (540,136.42) .. (540,132) -- cycle ;

%Straight Lines [id:da4250778022628171] 
\draw    (150,120) -- (268,120) ;
\draw [shift={(270,120)}, rotate = 180] [color={rgb, 255:red, 0; green, 0; blue, 0 }  ][line width=0.75]    (10.93,-3.29) .. controls (6.95,-1.4) and (3.31,-0.3) .. (0,0) .. controls (3.31,0.3) and (6.95,1.4) .. (10.93,3.29)   ;
%Straight Lines [id:da825761066474115] 
\draw    (540,120) -- (412,120) ;
\draw [shift={(410,120)}, rotate = 360] [color={rgb, 255:red, 0; green, 0; blue, 0 }  ][line width=0.75]    (10.93,-3.29) .. controls (6.95,-1.4) and (3.31,-0.3) .. (0,0) .. controls (3.31,0.3) and (6.95,1.4) .. (10.93,3.29)   ;
%Rounded Rect [id:dp45496863128326837] 
\draw  [fill={rgb, 255:red, 80; green, 227; blue, 194 }  ,fill opacity=1 ] (10,198) .. controls (10,193.58) and (13.58,190) .. (18,190) -- (142,190) .. controls (146.42,190) and (150,193.58) .. (150,198) -- (150,222) .. controls (150,226.42) and (146.42,230) .. (142,230) -- (18,230) .. controls (13.58,230) and (10,226.42) .. (10,222) -- cycle ;

%Straight Lines [id:da8144114572967778] 
\draw    (340,150) -- (340,170) -- (80,170) -- (80,188) ;
\draw [shift={(80,190)}, rotate = 270] [color={rgb, 255:red, 0; green, 0; blue, 0 }  ][line width=0.75]    (10.93,-3.29) .. controls (6.95,-1.4) and (3.31,-0.3) .. (0,0) .. controls (3.31,0.3) and (6.95,1.4) .. (10.93,3.29)   ;
%Rounded Rect [id:dp1763957119837839] 
\draw  [fill={rgb, 255:red, 155; green, 155; blue, 155 }  ,fill opacity=1 ] (270,196) .. controls (270,187.16) and (277.16,180) .. (286,180) -- (394,180) .. controls (402.84,180) and (410,187.16) .. (410,196) -- (410,224) .. controls (410,232.84) and (402.84,240) .. (394,240) -- (286,240) .. controls (277.16,240) and (270,232.84) .. (270,224) -- cycle ;
%Straight Lines [id:da3226973681574459] 
\draw    (270,210) -- (152,210) ;
\draw [shift={(150,210)}, rotate = 360] [color={rgb, 255:red, 0; green, 0; blue, 0 }  ][line width=0.75]    (10.93,-3.29) .. controls (6.95,-1.4) and (3.31,-0.3) .. (0,0) .. controls (3.31,0.3) and (6.95,1.4) .. (10.93,3.29)   ;
%Straight Lines [id:da5458300026642855] 
\draw    (610,140) -- (610,210) -- (412,210) ;
\draw [shift={(410,210)}, rotate = 360] [color={rgb, 255:red, 0; green, 0; blue, 0 }  ][line width=0.75]    (10.93,-3.29) .. controls (6.95,-1.4) and (3.31,-0.3) .. (0,0) .. controls (3.31,0.3) and (6.95,1.4) .. (10.93,3.29)   ;
%Rounded Rect [id:dp6737581508014832] 
\draw  [fill={rgb, 255:red, 155; green, 155; blue, 155 }  ,fill opacity=1 ] (270,294.5) .. controls (270,286.49) and (276.49,280) .. (284.5,280) -- (395.5,280) .. controls (403.51,280) and (410,286.49) .. (410,294.5) -- (410,325.5) .. controls (410,333.51) and (403.51,340) .. (395.5,340) -- (284.5,340) .. controls (276.49,340) and (270,333.51) .. (270,325.5) -- cycle ;
%Rounded Rect [id:dp4958062549820291] 
\draw  [fill={rgb, 255:red, 80; green, 227; blue, 194 }  ,fill opacity=1 ] (10,298) .. controls (10,293.58) and (13.58,290) .. (18,290) -- (142,290) .. controls (146.42,290) and (150,293.58) .. (150,298) -- (150,322) .. controls (150,326.42) and (146.42,330) .. (142,330) -- (18,330) .. controls (13.58,330) and (10,326.42) .. (10,322) -- cycle ;

%Straight Lines [id:da9731165295514408] 
\draw    (150,310) -- (268,310) ;
\draw [shift={(270,310)}, rotate = 180] [color={rgb, 255:red, 0; green, 0; blue, 0 }  ][line width=0.75]    (10.93,-3.29) .. controls (6.95,-1.4) and (3.31,-0.3) .. (0,0) .. controls (3.31,0.3) and (6.95,1.4) .. (10.93,3.29)   ;
%Rounded Rect [id:dp7966627240707042] 
\draw  [fill={rgb, 255:red, 245; green, 166; blue, 35 }  ,fill opacity=1 ] (540,298) .. controls (540,293.58) and (543.58,290) .. (548,290) -- (672,290) .. controls (676.42,290) and (680,293.58) .. (680,298) -- (680,322) .. controls (680,326.42) and (676.42,330) .. (672,330) -- (548,330) .. controls (543.58,330) and (540,326.42) .. (540,322) -- cycle ;

%Straight Lines [id:da5815062147793697] 
\draw    (410,310) -- (538,310) ;
\draw [shift={(540,310)}, rotate = 180] [color={rgb, 255:red, 0; green, 0; blue, 0 }  ][line width=0.75]    (10.93,-3.29) .. controls (6.95,-1.4) and (3.31,-0.3) .. (0,0) .. controls (3.31,0.3) and (6.95,1.4) .. (10.93,3.29)   ;

% Text Node
\draw (338.83,29.5) node   [align=left] {Meta-problem};
% Text Node
\draw (78.83,29.5) node   [align=left] {Input};
% Text Node
\draw (608.83,29.5) node   [align=left] {Output};
% Text Node
\draw (78.83,119.5) node   [align=left] {Parameters};
% Text Node
\draw (608.83,119.5) node   [align=left] {Key Clinical Output};
% Text Node
\draw (539,110.5) node [anchor=east] [inner sep=0.75pt]   [align=left] {define};
% Text Node
\draw (78.83,209.5) node   [align=left] {Relevant parameters};
% Text Node
\draw (88,180.5) node [anchor=west] [inner sep=0.75pt]   [align=left] {select};
% Text Node
\draw (180,209) node [anchor=south] [inner sep=0.75pt]   [align=left] {estimate};
% Text Node
\draw (151,110.5) node [anchor=west] [inner sep=0.75pt]   [align=left] {get ranges};
% Text Node
\draw (611,149.5) node [anchor=west] [inner sep=0.75pt]   [align=left] {measure};
% Text Node
\draw (78.83,309.5) node   [align=left] {Relevant parameters};
% Text Node
\draw (608.83,309.5) node   [align=left] {Key Clinical Output};
% Text Node
\draw (538,310) node [anchor=east] [inner sep=0.75pt]   [align=left] {\begin{minipage}[lt]{49.64pt}\setlength\topsep{0pt}
\begin{flushright}
confidently\\predict
\end{flushright}

\end{minipage}};
% Text Node
\draw (151,311) node [anchor=west] [inner sep=0.75pt]   [align=left] {estimate\\errors};
% Text Node
\draw  [fill={rgb, 255:red, 255; green, 255; blue, 255 }  ,fill opacity=1 ]  (286.5, 106.5) circle [x radius= 13.73, y radius= 13.73]   ;
\draw (286.5,106.5) node   [align=left] {1};
% Text Node
\draw (338.83,119.25) node   [align=left] {\begin{minipage}[lt]{49.64pt}\setlength\topsep{0pt}
\begin{center}
Sensitivity \\Analysis
\end{center}

\end{minipage}};
% Text Node
\draw  [fill={rgb, 255:red, 255; green, 255; blue, 255 }  ,fill opacity=1 ]  (287.5, 197.5) circle [x radius= 13.73, y radius= 13.73]   ;
\draw (287.5,197.5) node   [align=left] {2};
% Text Node
\draw (338.83,209.25) node   [align=left] {\begin{minipage}[lt]{38.76pt}\setlength\topsep{0pt}
\begin{center}
Inverse\\problem
\end{center}

\end{minipage}};
% Text Node
\draw  [fill={rgb, 255:red, 255; green, 255; blue, 255 }  ,fill opacity=1 ]  (286.5, 296.5) circle [x radius= 13.73, y radius= 13.73]   ;
\draw (286.5,296.5) node   [align=left] {3};
% Text Node
\draw (338.83,309.25) node   [align=left] {\begin{minipage}[lt]{63.24pt}\setlength\topsep{0pt}
\begin{center}
Uncertainty \\Quantification
\end{center}

\end{minipage}};

\end{tikzpicture}}
    \caption{
        \textbf{Schematic representation of the suggested pipeline for 0D model meta-analysis. An application of this pipeline is presented in section \ref{section:application}} 
    }
    \label{fig_pipeline}
\end{figure}
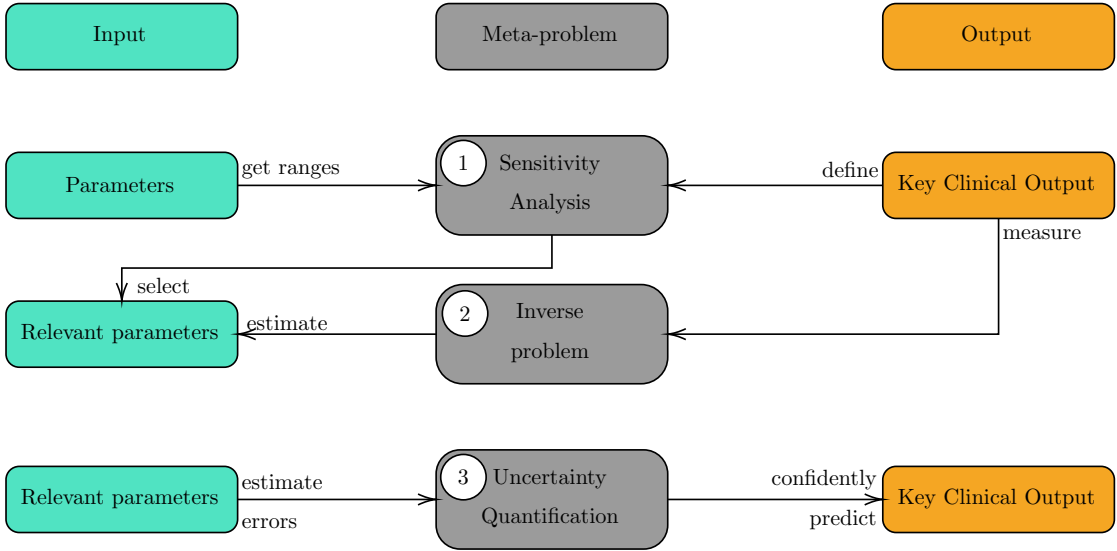

The main contributions of this work can be summarized as:
\begin{enumerate}
    \item Introducing a pipeline to build surrogate models for 0D cardiovascular models, which is beneficial for conducting sensitivity analysis, parameter estimation, and uncertainty quantification, as shown in figure~\ref{fig_pipeline}.\\
    \item Comparing the neural networks-based surrogates with models built with PCE and GP for the three cardiovascular models.\\
    \item Applying the full pipeline to one of the 0D models with neural networks.
    
\end{enumerate}

\section{0D cardiovascular models}\label{section:0D}

Zero-dimensional models, otherwise called lumped parameter models, are reduced order models that are solely time-dependent and are typically characterized by sets of (possibly non-linear) DAEs \cite{formaggia2010cardiovascular, shi2011review}. 

With regards to cardiovascular modelling, the key feature of the 0D models is their capability to simulate hemodynamics throughout the entire circulatory system, with a reasonable degree of accuracy. They bring the potential to offer clinicians a reliable tool for fast decision making \cite{golse2021predicting}.

Within the framework of 0D modeling, the hydraulic-electric analogy is employed, where blood flow rate ($Q$) corresponds to electrical current and blood pressure ($P$) to voltage. The detailed analogy between the different electrical components that are used in the three models and the pressure-flow relationship is presented in Table \ref{tab:pressure_flow_relationships}.

\begin{table}[!ht]
    \caption{Pressure-Flow relationships of the most commonly used electrical components. Resistance ($R$) accounts for the viscous loss effects of blood flow, while the non-linear (quadratic) resistance ($K$) captures non-linear effects caused by factors such as flow separations, stenoses, and other complex phenomena. Inertance ($L$) reflects the inertia of blood, representing the resistance to changes in flow. Compliance ($C$) describes the ability of a vessel to expand and store blood, while time-varying compliance accounts for the dynamic changes in compliance, typically replicating cardiac chamber dynamics \cite{Haghebaert2025}. The valve allows flow ($q$) only when the pressure difference ($p$) is positive, preventing backflow. Note that the illustrated valve model is an ideal representation.}
    \centering
    \renewcommand{\arraystretch}{1.5} % Adjust row height
    \begin{tabular}{|c|c|}
        \hline
        \textbf{Component} & \textbf{Pressure-Flow Relationship} \\
        \hline
        \begin{minipage}{0.4\textwidth}
            \centering
            \vspace{0.2cm} % Add vertical space
            Resistance \\
            \includegraphics[width=0.4\textwidth]{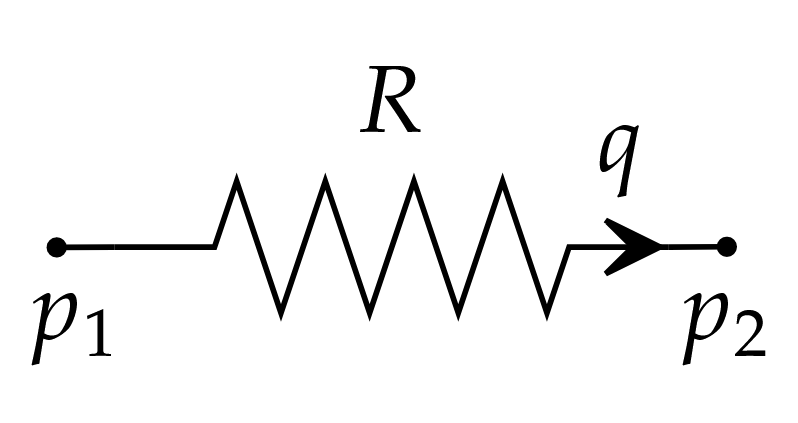}
        \end{minipage} &
        $p_{1} - p_{2} = Rq$ \\
        \hline
        \begin{minipage}{0.4\textwidth}
            \centering
            \vspace{0.2cm} % Add vertical space
            Non-linear Resistance \\
            \includegraphics[width=0.4\textwidth]{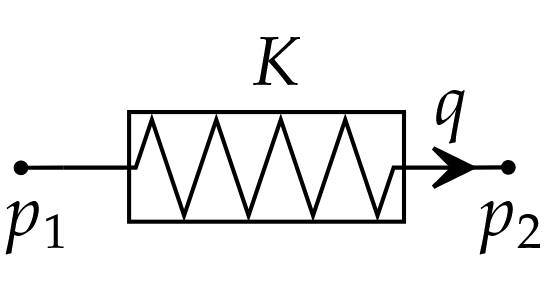}
        \end{minipage} &
        $p_{1} - p_{2} = Kq|q|$ \\
        \hline
        \begin{minipage}{0.4\textwidth}
            \centering
            \vspace{0.2cm} % Add vertical space
            Inertance \\
            \includegraphics[width=0.4\textwidth]{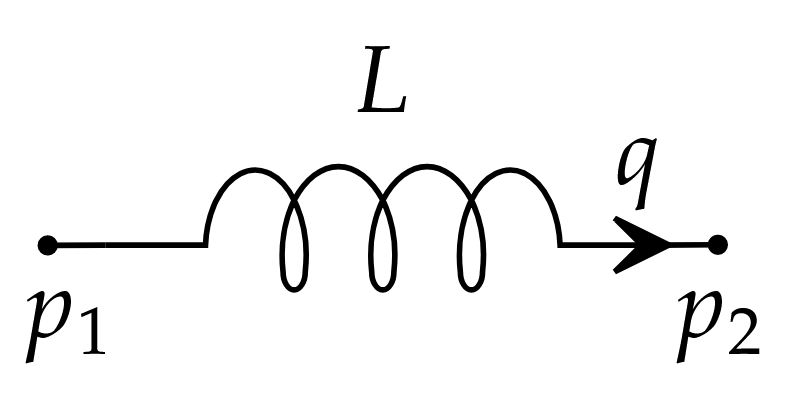}
        \end{minipage} &
        $p_{1} - p_{2} = L \dot{q}$ \\
        \hline
        \begin{minipage}{0.4\textwidth}
            \centering
            \vspace{0.2cm} % Add vertical space
            Compliance \\
            \includegraphics[width=0.4\textwidth]{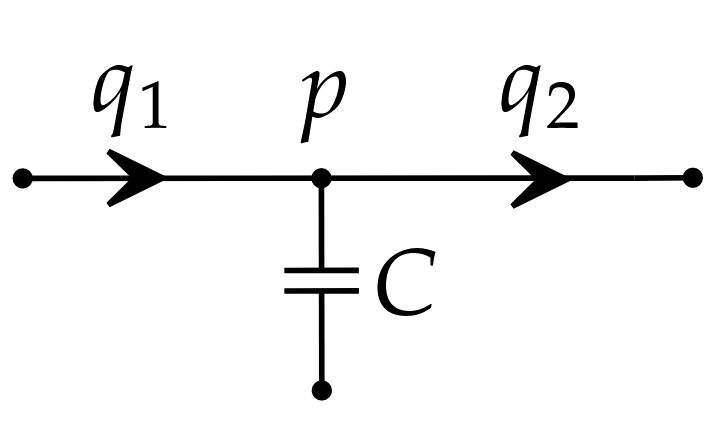}
        \end{minipage} &
        $q_{1} - q_{2} = C \dot{p}$ \\
        \hline
        \begin{minipage}{0.4\textwidth}
            \centering
            \vspace{0.2cm} % Add vertical space
            Time-varying Compliance \\
            \includegraphics[width=0.4\textwidth]{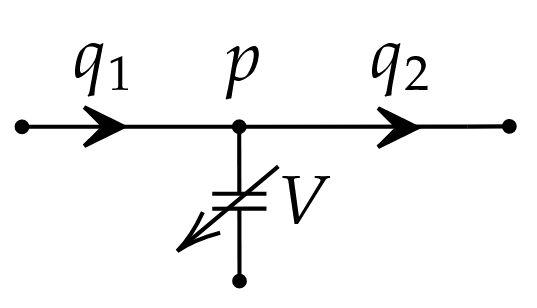}
        \end{minipage} &
        $q_{1} - q_{2} = \dot{v}$ \\
        \hline
        \begin{minipage}{0.4\textwidth}
            \centering
            \vspace{0.2cm} % Add vertical space
            Valve \\
            \includegraphics[width=0.4\textwidth]{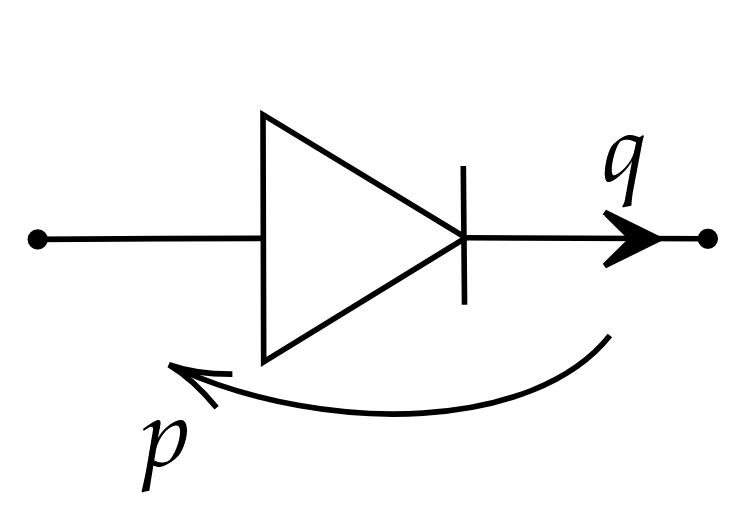}
        \end{minipage} &
        \makecell{ 
            $p = 0$ if $q > 0$, \\
            $q = 0$ if $p < 0$
        } \\
        \hline
    \end{tabular}\label{tab:pressure_flow_relationships}
\end{table}

In Table \ref{tab:pressure_flow_relationships}, the non-linear resistance is a concept not often observed in 0D models, yet, in our case, it is included in the first two models. It takes into account the nonlinear effects of convective energy loss due to possible flow separation, vortices, and turbulence, mainly occurring in larger vessels \cite{migliavacca2001modeling}. Additionally, Model 1 contains a pressure-dependent resistance, which describes an increase in vessel radii due to increased pressure, thus decreasing vascular resistance \cite{audebert2017partial}. Lastly, for a more dynamic representation of the valve component, the addition of a non-linear resistance and an inductance is considered, as proposed by \cite{mynard2012simple}, to replicate the leaflet motion. More detailed expressions for valve dynamics have also been introduced, for instance, considering regurgitation and stenosis \cite{pant2018lumped, pant2016data}.

From a numerical perspective, a lumped parameter model represents a system of differential-algebraic equations (DAE) given by:

\begin{equation}
\left\{
\begin{aligned}
    &\frac{dy}{dt} = b(y, z, t), & t &\in (0, T], \\
    &        G(y, z, t) = 0, &
\end{aligned}
\right .
\label{eq:DAE}
\end{equation}

where the initial condition vector is \( y|_{t=t_0} = y_0 \). In Eq.\eqref{eq:DAE}, \(y\) denotes the vector of state variables, \(z\) encompasses the remaining network variables, and \(G\) states the algebraic equations. The numerical methods to solve DAEs usually involve a combination of techniques for handling the differential and algebraic components of the system. The Runge-Kutta iterative method is often employed.

While 0D models typically offer rapid computation and reliable solutions, they frequently suffer from over-parameterization. This means that they include parameters not directly measured in clinical settings, necessitating tuning or assumptions based on analogous studies.

\subsection{Model 1: predicting portal vein pressure after liver surgery}
Model 1 aims at predicting hemodynamic changes after partial hepatectomy, whereby a significant portion of liver tissue is removed \cite{audebert2017partial, golse2021predicting}. Removal of tissue leads to an increase of resistance for mesenteric blood flow, which in turn will lead to an increased portal pressure. Current guidelines consider pre-operative portal hypertension (porto-caval pressure gradient $\geq$ 10 mmHg) \cite{liver2012easl} before surgery as a risk for the potential failure of the remaining liver tissue, however, \citep{allard2013posthepatectomy} instead demonstrated post-operative portal hypertension (portal pressure $\geq$ 21 mmHg) to be a better predictor of post-operative liver failure. Hence, model-based prediction of portal pressure may alleviate risks of post-operative complications. 

\begin{figure}[H]
    \centering
    \resizebox{0.8\textwidth}{!}{\input{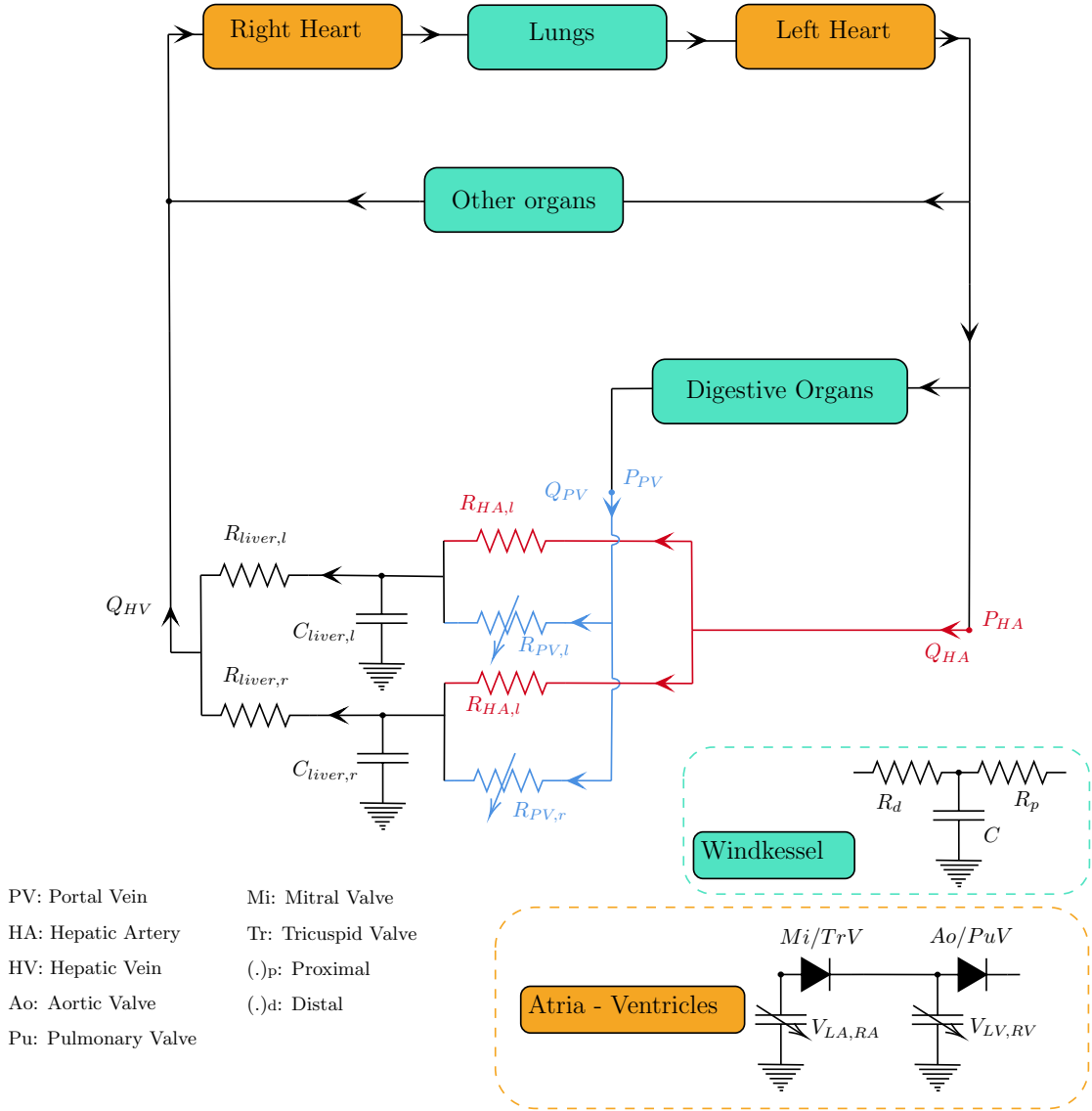}}
    \caption{
        \textbf{Schematic for Model 1, for prediction of hemodynamic change after liver resection \cite{audebert2017partial, golse2021predicting} It is composed of the left and right parts of the heart, pulmonary circulation, digestive organs, two hemi-livers and remaining circulation. } 
    }
    \label{fig_LumpedFlowDiagram}
\end{figure}

Since approximately 25\% of the total cardiac output passes the liver, changes in liver resistances by removal of tissue will likely alter the total hemodynamics of the system. Therefore, the model comprises the total circulation (simplified): left and right heart chambers and associated valves, pulmonary circulation, digestive organs, two hemi-livers and remaining circulation (Figure \ref{fig_LumpedFlowDiagram}). The model is personalized with pre-operative data, and resection is simulated by increasing resistance proportional to the fraction of tissue removed \cite{audebert2017partial, golse2021predicting}.

\subsection{Model 2: pulmonary arterial hypertension intervention planning}
The second model is targeted for simulating the hemodynamic conditions under Pulmonary Arterial Hypertension (PAH) \cite{farber2004pulmonary}. PAH is a rare disease, which can be idiopathic or related to secondary diseases, and it is described by an abnormal increase in the pulmonary arterial pressure. It can lead to right ventricular overloading, progressively inducing a right ventricular failure, and if left untreated, it can be life-threatening \cite{valdeolmillos2023thirty}. This zero-dimensional model, originally proposed in \cite{pant2022multiscale}, serves as a predictive tool for intervention planning. When drug therapy proves ineffective, procedures such as the creation of a pulmonary-to-systemic shunt, commonly known as Potts shunt \cite{blanc2004potts}, provides an alternative to lung transplantation for a specific group of patients.

In terms of its technical characteristics, the main distinction between this model and the first one outlined in this study lies in its utilization of the single fiber model  \cite{arts1991relation, bovendeerd2006dependence} to simulate cardiac electromechanics, as opposed to the time-varying elastance model employed in the former. In addition to that, the upper and lower body circulations are represented more finely, while the shunt region, including the ascending and descending aorta, the aortic arch, and the right and left pulmonary arteries, are replicated in detail, in contrast to the first focusing on the liver circulation. This closed-loop circulation model can also be coupled with a 3D model for performing CFD simulations in the shunt region as shown in \cite{pant2022multiscale}, in a similar way to has already been done for several other applications within the context of cardiovascular modelling \cite{marsden2015multiscale, lagana2005multiscale, migliavacca2006multiscale}. The 3D flow domain is also initially required to tune this 0D model as described by \cite{pant2014multiscale, pant2014methodological}.

The model can capture the hemodynamics of the hypertensive patients before and after the shunt integration between the left pulmonary artery and the descending aorta. Personalization based on preoperative patient data is a non-trivial task \cite{pant2022multiscale}. Several clinically relevant key parameters can then be extracted from the model to assess the condition of the patient. Some of these are pressures, systemic and pulmonary flow rates, working volumes, and ejection fractions for assessing cardiac function and post shunt integration. 

\begin{figure}[H]
    \centering
    \resizebox{0.8\textwidth}{!}{\input{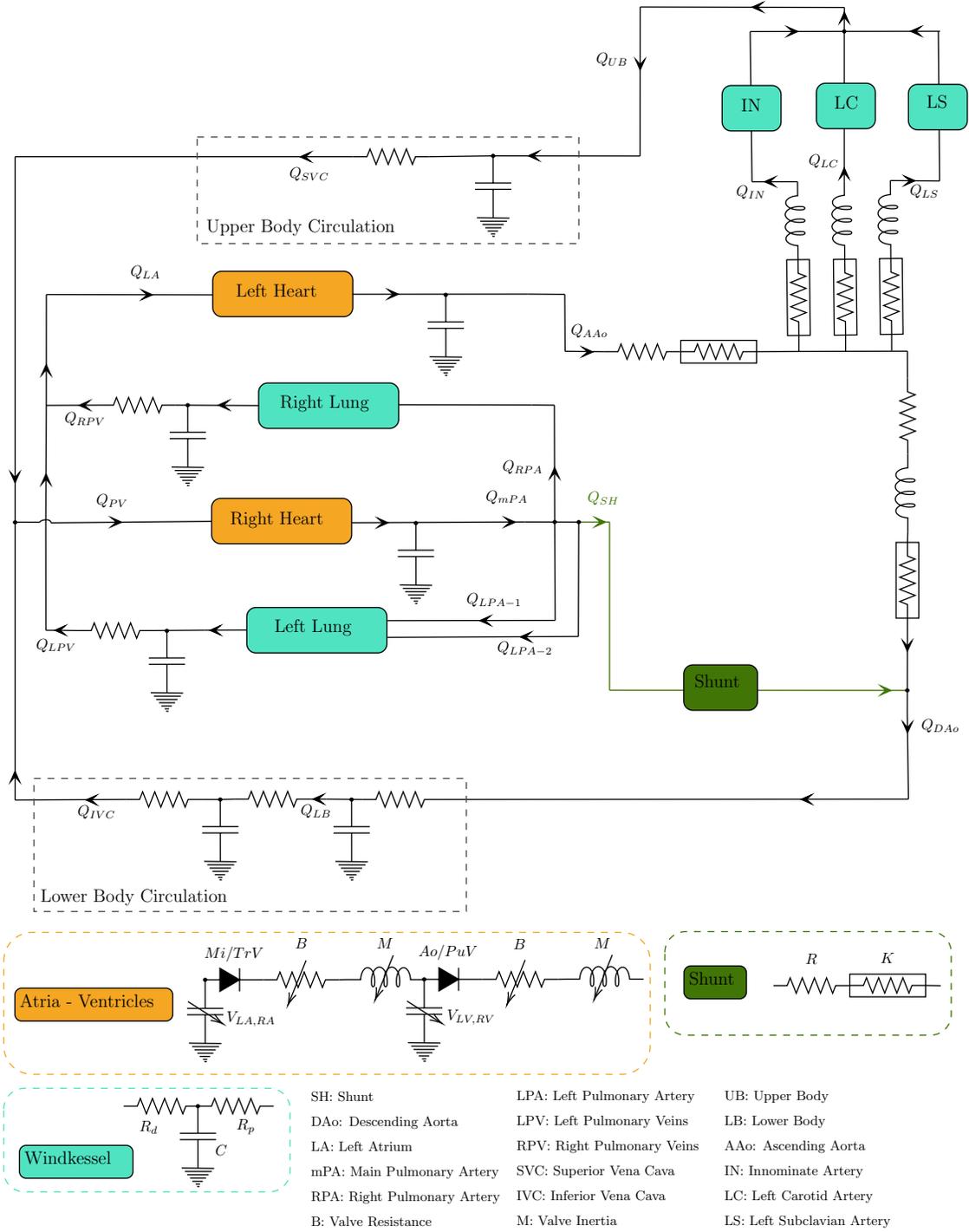}}
    \caption{
        Schematic diagram of the PAH model with integrated conduit shunt. $Q$ represents the volumetric flow-rate. The non-linear resistances are extracted after tuning with the 3D model of the shunt region. The heart and circulatory system are studied in more details compared to the first model.
    }
    \label{fig:PAHmodel}
\end{figure}

\subsection{Model 3: assessment of perfusion through fluorescence patterns}
 The third model focuses on the time variation of a contrast agent in a certain region of interest for medical imaging applications. In addition to a 0D hemodynamic model (figure \ref{fig_LumpedFlowDiagram}, table \ref{tab:pressure_flow_relationships}), the transport of a substance is evaluated. The substance is injected intravenously (IV injection) and is eliminated at a constant rate (\textit{e.g.} through filtration in the kidneys or metabolization in the liver). It is transported in the blood through the main circulation and reaches the organ of interest, where the fluorescence intensity is measured. A vascular anomaly compartment is inserted before the organ of interest to model an abnormal reduction of blood flow going to the organ of interest and that results in a reduction of the organ perfusion.

\begin{figure}[H]
    \centering
    \normalsize
        \resizebox{0.8\textwidth}{!}{\input{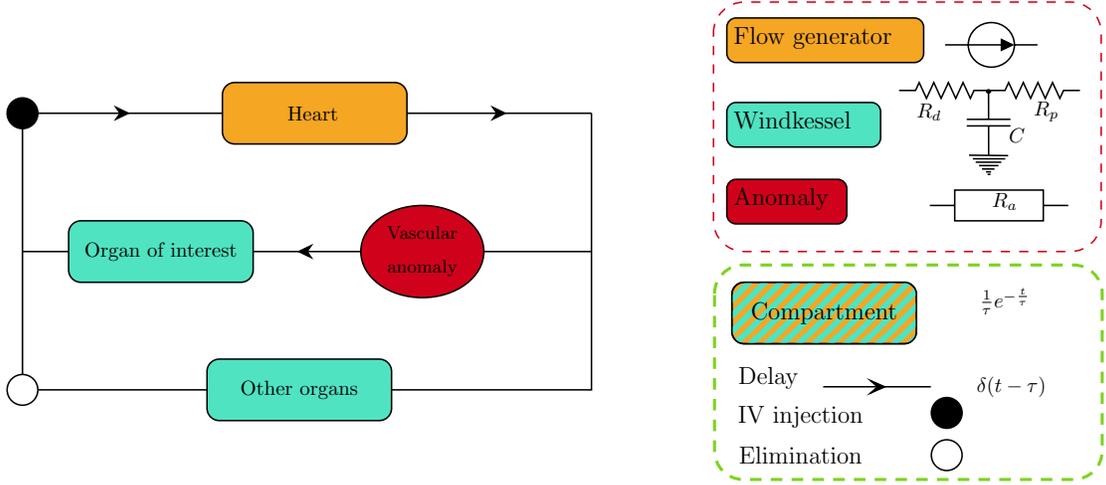}}
    \caption{Schematic diagram of the hemodynamic and transport model for assessing perfusion in the organ of interest. The top red dashed box is the key for the hemodynamics model and the bottom green box is the key to the transport model. Here the heart is modelled by a simple flow generator (see above for more precise heart models). Hemodynamics is solved before transport. The $h(t)$ function represents the transformation undergone by the compartment input signal through the convolution relation $c_\mathrm{out} = h * c_\mathrm{in}$. Parameters $\tau$ depend on the organs and the hemodynamics properties.}
    \label{fig:compartments}
\end{figure}

The advantage of the model is that the vascular anomaly is represented through a geometric parameter that aims to replicate the behaviors of the contrast agent in different vascular configurations. The objective is for the metamodel to serve as a surrogate for the inverse problem: retrieving the value of the anomaly parameter with the sole information of the time variation of the contrast agent concentration in the organ of interest. The specificity of this study lies in the fact that the input of the metamodel is a time series (in contrast to the other models, where there is a discrete set of inputs).

\section{Methodology}\label{section:nn_method}

\subsection{Surrogate modelling}

\subsubsection{Feed-forward neural networks}

Feed-forward neural networks (FFNN) or sometimes referred to as fully connected neural networks, are a class of approximation functions that are commonly applied in supervised learning tasks. They are widely used due to their strong approximation capabilities, where according to the universal approximation theorem, neural networks can approximate any continuous function with accuracy depending on the optimization. The FFNN approximation has the form of

\begin{equation}
     \mathbf{z}^i = \sigma^i(\mathbf{W}^i \mathbf{z}^{i-1} + \mathbf{b}^i),\ \ i=1,...,L
 \end{equation}

 where $\mathbf{z}^0$ and $\mathbf{z}^L$ are the input and output of the function, respectively. While, $\mathbf{W}^i$ and $\mathbf{b}^i$ are the parameters of each layer, known as the weights and biases, respectively. $\sigma^i$ is a nonlinear function used to add nonlinearity in the representation and is called the activation function. $L$ is the number of layers in the network.

Given $N$ data points in the form of ${(\mathbf{X}_i, \mathbf{Y}_i)}_{i=1:N}$, neural networks can be used to perform a supervised learning task by minimizing a chosen loss function, a mean squared error for example, with respect to the neural network parameters.

\begin{equation}
    \mathcal{L} = \sum_{i=1}^N |\mathbf{\hat{Y}}(\mathbf{X}_i; \mathbf{\theta})-\mathbf{Y}_i)|^2,
\end{equation}

where $\mathbf{\hat{Y}}$ is a neural network approximation and $\mathbf{\theta}$ a vector comprising the set of function parameters (weights and biases).

FFNN is well suited for approximating tabular or scalar data, however, that is not the case for time series data. That is because the sequential nature of the data is ignored in the architecture and it has no memory mechanism. Moreover, the inputs or outputs of the network need to have a fixed size.

Recurrent neural network architectures \cite{rnn} are designed to treat time-series data where the sequential nature of the data is the key in the architecture design. Moreover, some variants such as Long-Short Term Memory (LSTM) network \cite{lstm} have long term dependencies and address exploding/vanishing gradients issue.

In this work, the neural networks are implemented in the Tensorflow library \cite{tensorflow2015-whitepaper}.

\subsubsection{Polynomial chaos expansion}
Polynomial Chaos Expansion (PCE) was originally designed by \citep{wiener1938homogeneous} to represent random variables in the form of a set of polynomials and is often utilized in global sensitivity analysis \cite{sudret2008global}. The choice of polynomial corresponds to the probability distribution of its respective model input parameter. It is possible to approximate the response $Y$ of a system $f$ to a random input $\mathbf{X}$ using the aforementioned polynomials,

\begin{equation}
    Y = f(\mathbf{X}) \approx f_{PCE}(\mathbf{X}) = \sum_{j=0}^{P-1} \alpha_j \Psi_j(\mathbf{X})
    \label{eq:PCE}
\end{equation}

with $\Psi$ and $\alpha$ the set of orthonormal polynomials and their coefficients. Coefficients can be obtained via several regression methods, e.g. least-squares, ridge or least angle regression. The number of samples $N$ required to find a unique solution for \ref{eq:PCE} is dependent on the number of input parameters $D$ and the maximal order polynomial order $P$, given by
\begin{equation}
    N = \binom{D+P}{D}.
\end{equation}
However, since interaction terms are likely to contribute significantly less, several truncation schemes exist \cite{blatman2008sparse} to reduce the number of samples required. An important novel addition is the so-called \textit{adaptive sparse} PCE, which iteratively adds and removes polynomial terms based on their contribution to the surrogate model's output of interest, thereby simultaneously reducing the number of input samples required as well as reducing the chance of overfitting \cite{blatman2011adaptive}. In this case, model quality is given by the jack-knife error, whereby the quality is determined by sequentially removing a single data point and determining the model prediction error with respect to that point. 

A benefit of PCE is that the Sobol sensitivity indices can be derived analytically from the PCE coefficients \cite{sudret2008global}, informing the user which input parameters might be fixed parameters. 

In this work, we use the Python implementation of \textit{OpenTURNS} \cite{OpenTURNS}, which utilises the adaptive PCE.

\subsubsection{Gaussian processes}
Gaussian Processes (GPs) are a powerful statistical tool for modeling and understanding complex datasets \cite{marrel2024probabilistic}. 
At their core, GPs are a type of Bayesian inference method that provides a flexible approach to regression and classification problems. They excel in situations where the underlying structure of the data is unknown or difficult to specify with traditional parametric models.

A GP defines a distribution over functions where any finite collection of these function values has a joint Gaussian distribution. This characteristic makes GPs incredibly versatile and capable of modeling a wide range of functions. The heart of a GP model lies in its kernel or covariance function, which encodes assumptions about the function being modeled, such as smoothness or periodicity.
By choosing and tweaking this kernel function, one can control the behavior and adaptiveness of the GP.
A GP is formally defined as a collection of random variables, any finite number of which have a joint Gaussian distribution. This can be succinctly represented as:
\begin{equation*}
    f(x) \sim \mathcal{GP}(m(x), k(x, x'))
\end{equation*}
where $f(x)$ is the function to be modeled, $ m(x) $ is the mean function, which is usually assumed to be zero because the GP is often more concerned with learning deviations from a mean rather than the mean itself, and $ k(x, x')$ is the kernel or covariance function, which describes the relationship between any two points in the input space.
The kernel function is a critical component of the GP framework, as it encodes the assumptions about the function to be modeled. 
A common choice for the kernel is the Squared Exponential (SE) or Gaussian kernel:
\begin{equation}
    k(x, x') = \sigma^2 \exp{\left(-\frac{\|x - x'\|^2}{2l^2}\right)}
\end{equation}
where $\sigma^2$ is the variance parameter, controlling the overall variance of the process,  $l$ is the length scale parameter, governing the smoothness of the function, and $\|x - x'\| $ is the Euclidean distance between points $x$ and $x'$.

Given a set of training data points $(X, y)$ and a new test input $x^*$, the GP prediction for $x^*$ is characterized by the predictive distribution:
\begin{equation*}
    f(x^*)|X, y, x^* \sim \mathcal{N}(\mu^*, \sigma^{*2})
\end{equation*}
with:
\begin{align*}
    \text{mean} \quad & \mu^*  = k(x^*, X)[k(X, X) + \sigma_n^2I]^{-1}y \\
    \text{variance} \quad & \sigma^{*2} = k(x^*, x^*) - k(x^*, X)[k(X, X) + \sigma_n^2I]^{-1}k(X, x^*).
\end{align*}
Here $ k(X, X)$ is the covariance matrix computed from the training inputs $ X $, $ \sigma_n^2 $ is the noise term in the observations, and $ I $ is the identity matrix. 
This setup allows GPs not only to predict the most likely value $\mu^*$ but also to provide an estimate of the uncertainty of this prediction $ \sigma^{*2}$.

In this work, we employed the Python library scikit-learn \cite{scikit-learn} to implement our Gaussian Process model. 
%GPs offer a powerful and adaptable approach for modeling complex data. By judiciously choosing the kernel function and efficiently managing computational challenges, one can leverage GPs to glean insightful predictions and understandings from a wide array of data types and distribution patterns.

\subsection{Data preparation and hyperparameter selection}
In this work, synthetic data from each model are the basis to train the metamodels. The generated data are split into training and testing sets. For all models, the data are generated by implementing the sampling method described by Saltelli's algorithm \cite{saltelli2004sensitivity}, which is a quasi-Monte Carlo approach, mainly used to calculate the Sobol sensitivity indices \cite{sobol2001global}. The open-source Python library SALib \cite{herman2017salib} was implemented. The maximum computational cost for this approach is $N_s = (2d+2)N$, where $N_s$ represents the number of simulations or model evaluations, $d$ signifies the model's input variables, and $N$ denotes the sample size.

For the architecture of the neural networks, three hidden layers are chosen with 64, 32, and 32 neurons, consecutively. The Relu activation function is applied for the hidden layers. The Adam minimization algorithm is chosen to minimize the loss function; moreover, early stopping is implemented to avoid overfitting. In addition to that, the whole dataset is normalized before the training process.

For the polynomial chaos expansion, an adaptive algorithm iteratively removes polynomial terms $\alpha_j \Psi_j$ with a low relative contribution, i.e. where the coefficient $\alpha_j$ is below a certain threshold

\begin{equation}
    \mathcal{A}_{rem} = \{ \|\alpha_j\| \leq \epsilon \max_{\alpha \in \mathcal{A}_{tot}} \boldsymbol{(\alpha)} \},
\end{equation}
with $\mathcal{A}_{tot}$ the total set of polynomials and $\mathcal{A}_{rem}$ the set of polynomials to be removed. The error term $\epsilon$ set to $10^{-4}$.

For GPs methodology, for all models, the first step involves normalizing the training input data using a standard scaling technique.
This ensures that each feature contributes equally to the model. Similarly, the output data for each output dimension is scaled to standardize its range. The GPs are created using the Radial Basis Function (RBF) kernel, renowned for its flexibility and adeptness at managing non-linear relationships. This model undergoes training with the normalized data, effectively learning the inherent patterns and connections.
For predictions on the testing set, the inputs are scaled using the same normalization applied to the training data, and the predictions are adjusted back to their original scale. This guarantees that the output maintains the same format and range as the original dataset. Utilizing scikit-learn's comprehensive Gaussian Processes implementation, this methodology promotes efficient training and accurate predictions.
\section{Metamodel building and assessment}\label{section:results}

The same training and testing datasets are used for neural networks, PCE, and GP, to ensure fair comparison. In all cases, 20\% of the data is kept for testing, and the $Q^2$ metric (Eq. \ref{eq:q2}) is measured using different training set sizes (same datasets for all methods) along with the maximum error. This is done to see the effect of the training data size on the metamodels' performance.

\begin{equation}\label{eq:q2}
    Q^2 = 1 - \frac{\sum_{i=1}^{n} \left( \mathbf{Y}_i - \mathbf{\hat{Y}} (\mathbf{X}_i)\right)^2}{\sum_{i=1}^{n} \left( \mathbf{Y}_i - \mathbf{\bar{Y}} \right)^2},
\end{equation}

where $\mathbf{\bar{Y}}$ is the mean of the true values $\mathbf{Y}_i$, and $n$ is the testing dataset size.

\subsection{Model 1}
\begin{figure}[H]
    \centering    \includegraphics[width=0.9\columnwidth]{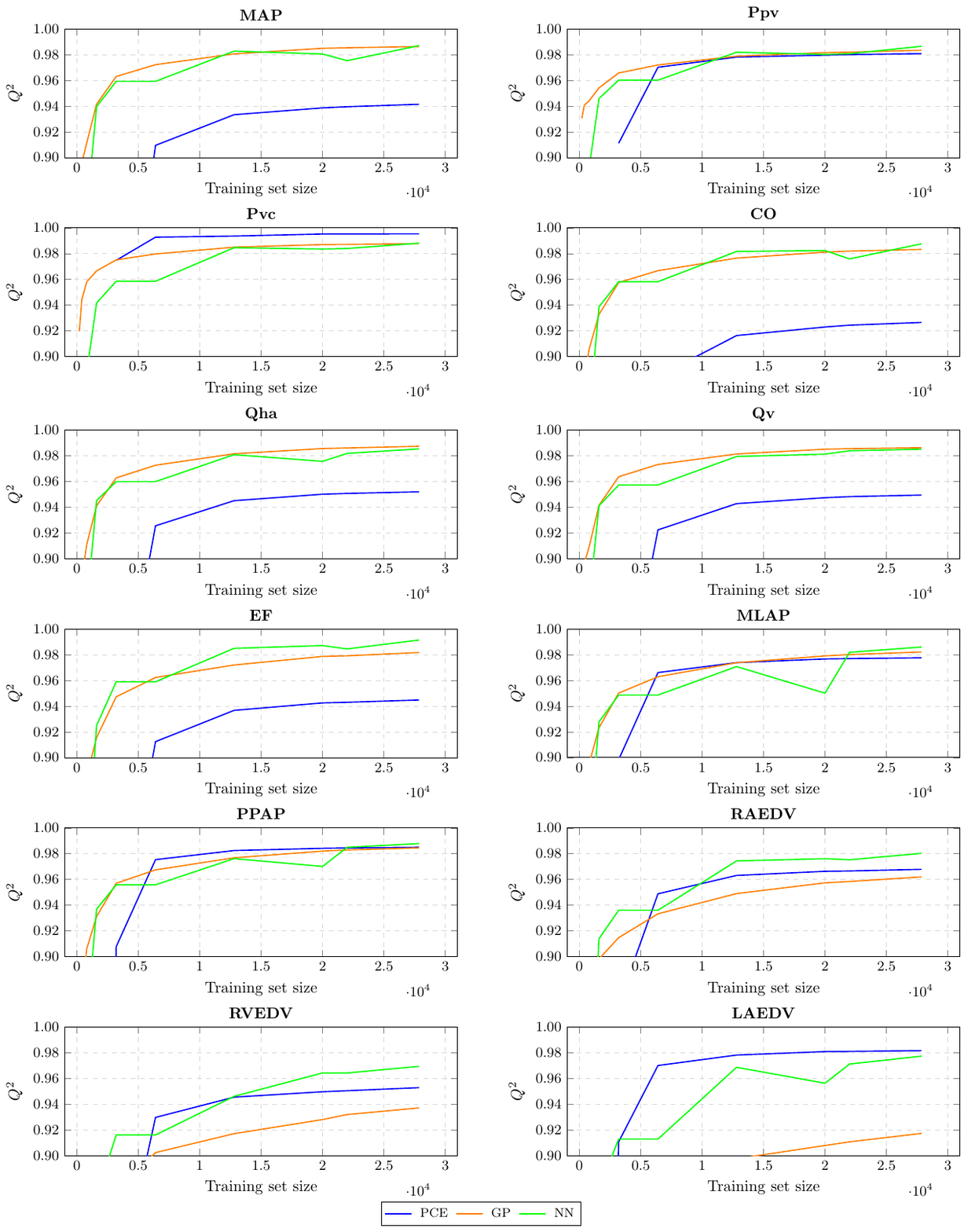}  %png. or .pdf for latex
    \caption{$Q^2$ metric using the testing dataset for the first model.}
    \label{fig:q2_roel}
\end{figure}

% After profile likelihood analysis \cite{vanlier2012integrated}, we determined elastances of the four heart chambers are not uniquely identifiable if there is no information on both the pressure and volume. Hence, we only estimate the left ventricular elastances while scaling the remaining chambers, similar to \cite{golse2021predicting}. 
Table \ref{tab:InputOutputModel1} presents the different types of model inputs as well as model outputs. Note that some parameters which are commonly considered outputs (e.g., Cardiac Output, Mean Aortic Pressure) are used as input to the model. This is due to the workflow whereby several resistances are directly derived from measured data, as explained in \cite{golse2021predicting}. The $Q^2$-metric for all outputs is plotted in figure~\ref{fig:q2_roel}.

\begin{table}[H]
\centering
\caption{Input and output parameters for Model 1}
\label{tab:InputOutputModel1}
\begin{tabular}{|l|l|}
\hline
\multicolumn{1}{|c|}{\textbf{Inputs}} & \textbf{Outputs} \\ \hline
Active Right Atrial Elastance ($Ea_{RA}$) & Mean Aortic Pressure, achieved ($MAP$)   \\ \hline
Baseline Right Atrial Elastance ($Eb_{RA}$) & Portal Vein Pressure, achieved ($P_{PV}$)   \\ \hline
Active Right Ventricular Elastance ($Ea_{RV}$) & Vena Cava Pressure, achieved ($P_{VC}$)  \\ \hline
Baseline Right Atrial Elastance ($Eb_{RV}$)  & Cardiac Output, achieved ($CO$) \\ \hline
Active Left Atrial Elastance ($Ea_{LA}$)  & Hepatic Artery Flow ($Q_{HA}$)  \\ \hline
Baseline Left Atrial Elastance ($Eb_{LA}$) & Portal Vein Flow ($Q_{PV}$) 
 \\ \hline
Active Left Ventricular Elastance ($Ea_{LV}$)  & Left Ventricular Ejection Fraction ($EF$)      \\ \hline
Baseline Left Ventricular Elastance ($Eb_{LV}$)   & Mean Left Atrial Pressure ($MLAP$)           \\ \hline
Vena Cava Compliance ($C_{vc}$) &    Peak Pulmonary Artery Pressure ($PPAP$)               \\ \hline
Mean Aortic Pressure, target ($MAP_m$) & Right Atrial End Diastolic Volume ($RAEDV$) \\ \hline
Portal Vein Pressure, target ($P_{PV,m}$) & Right Ventricular End Diastolic Volume ($RVEDV$) \\ \hline
Vena Cava Pressure, target ($P_{VC,m}$) & Left Atrial End Diastolic Volume ($LAEDV$) \\ \hline
Cardiac Output, target ($COm$) & \\ \hline
Hepatic Artery Flow, target ($Q_{HA,m}$) & \\ \hline
Portal Vein Flow, target ($Q_{PV,m}$) & \\ \hline
Left Liver Volume ($M_L$) & \\ \hline
Right Liver Volume ($M_R$) & \\ \hline
Heart Rate ($HR$) & \\ \hline
\end{tabular}
\end{table}

\begin{figure}[H]
    \centering    \includegraphics[width=0.9\columnwidth]{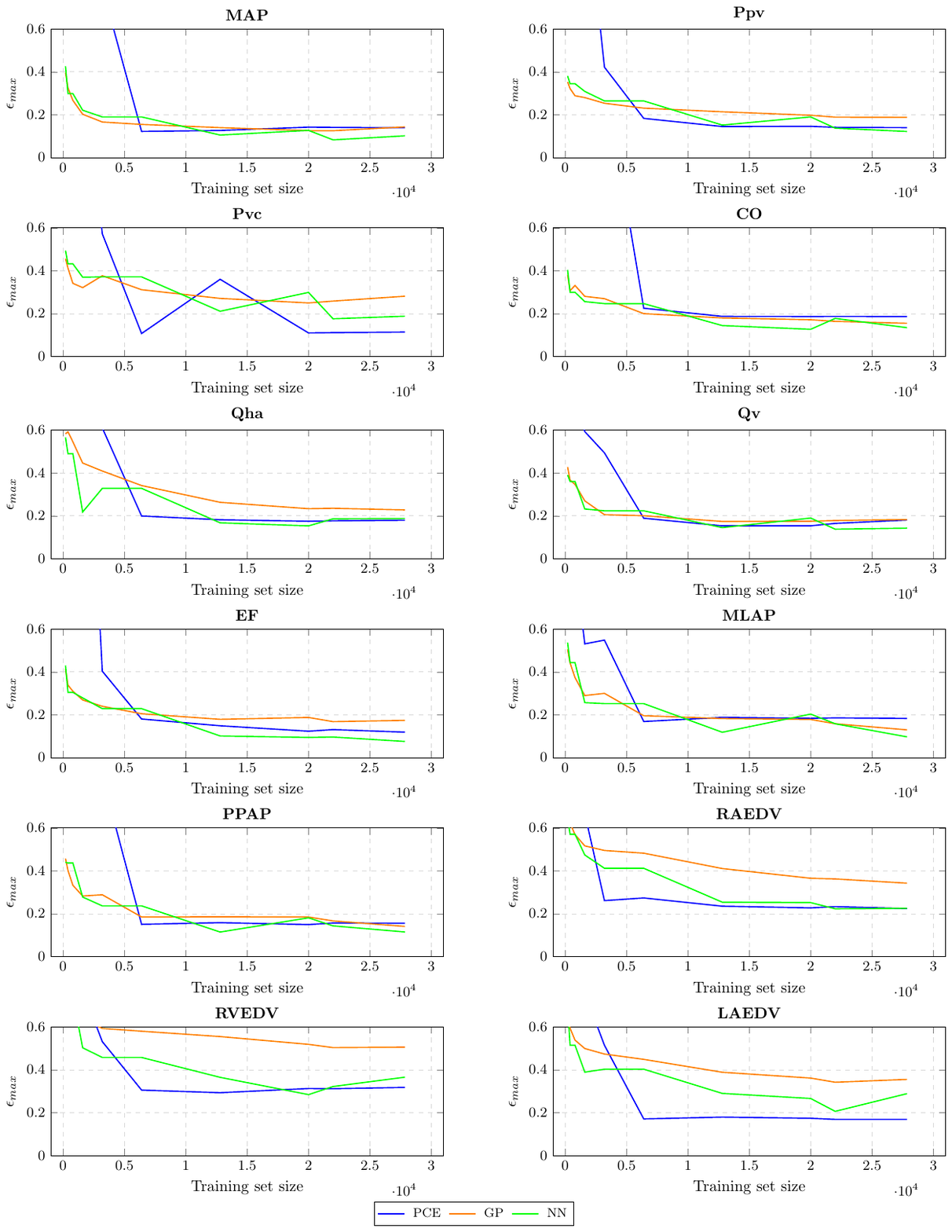}  %.png or .pdf for latex
    \caption{Maximum relative error using the testing dataset for the first model.}
    \label{fig:max_err_roel}
\end{figure}

As evident from figure~\ref{fig:q2_roel}, that PCE has poor performance in cases of low training datasets below 5000 points. For large training sets, NN and GP have better performance with neural networks showing some superiority in most of the outputs. To further analyze the results, another metric is shown in figure~~\ref{fig:max_err_roel}, which is the maximum relative error in the testing set.

From Figure~~\ref{fig:max_err_roel}, it can be noticed that GP has poor performance in the max relative error metric, which reaches above 50\% in the case of RVEDV output. This indicates that GP performs well on average; however, the standard deviation is high. On the contrary, we can see that PCE has a lower max relative error but a worse average error. Finally we can observe that NN has acceptable performance both in terms of the average and standard deviation error.  

\subsection{Model 2}\label{sec:model_2}

For this model, we have selected the ten most influential and clinically relevant inputs to the outputs based on the Morris Sensitivity Method \cite{morris1991factorial}. The Morris method performs a qualitative sensitivity analysis primarily suited for models with numerous independent input parameters. Given that the current model initially consists of over 40 independent inputs, many of which have minimal impact on the outcomes, the Morris method offers a cost-effective means to qualitatively rank them, in contrast to the Sobol indices method. With regard to the model outputs, the selection is based on the expertise of clinicians. The detailed selection of the model inputs and outputs for the metamodels' structure is presented in Table \ref{tab:model2inputsoutputs}

\begin{table}[htbp]
    \centering
    \caption{Input and output parameters for Model 2}
    \begin{tabular}{|c|c|}
    \hline
    \textbf{Inputs} & \textbf{Outputs} \\
    \hline
    Pulmonary Vascular Resistance ($PVR$) & Right Ventricular Stroke Volume ($RVSV$) \\
    \hline
    Systemic Vascular Resistance ($SVR$) & Left Ventricular Stroke Volume ($LVSV$) \\
    \hline
    Shunt Diameter ($d_{\text{shunt}}$) & Aortic mean pressure ($P_{{ao}_{mean}})$ \\
    \hline
    Coefficient for shunt's non-linear resistance ($k_{\text{shunt}}$) & Pulmonary Artery mean pressure ($P_{{pa}_{mean}})$ \\
    \hline
    Cardiac contractility and max fiber stress ($T_c$) & Mean Pressure Gradient across the shunt ($\Delta P_{{\text{shunt}_{\text{mean}}}}$) \\
    \hline
    Right Ventricular volume at zero transmural pressure ($V_{RV_{0}}$)  & Pulmonary Artery diastolic pressure ($P_{{pa}_{dias}})$ \\
    \hline
    Right Ventricular wall volume ($V_{w_{RV}}$)  & Pulmonary Artery systolic pressure ($P_{{pa}_{sys}})$ \\
    \hline
    Right Atrial wall volume ($V_{w_{RA}}$)  & Aortic systolic pressure ($P_{{ao}_{sys}})$ \\
    \hline
    Ratio between proximal and distal resistance ($r=Rp:Rd$)  & Aortic diastolic pressure ($P_{{ao}_{dias}})$ \\
    \hline
    Right Ventricular activation time ($t_{RV}$) & Pulmonary to Systemic flow ratio ($Qp:Qs$) \\
    \hline
    \end{tabular}
    \label{tab:model2inputsoutputs}
\end{table}

\begin{figure}[H]
    \centering    \includegraphics[width=0.9\columnwidth]{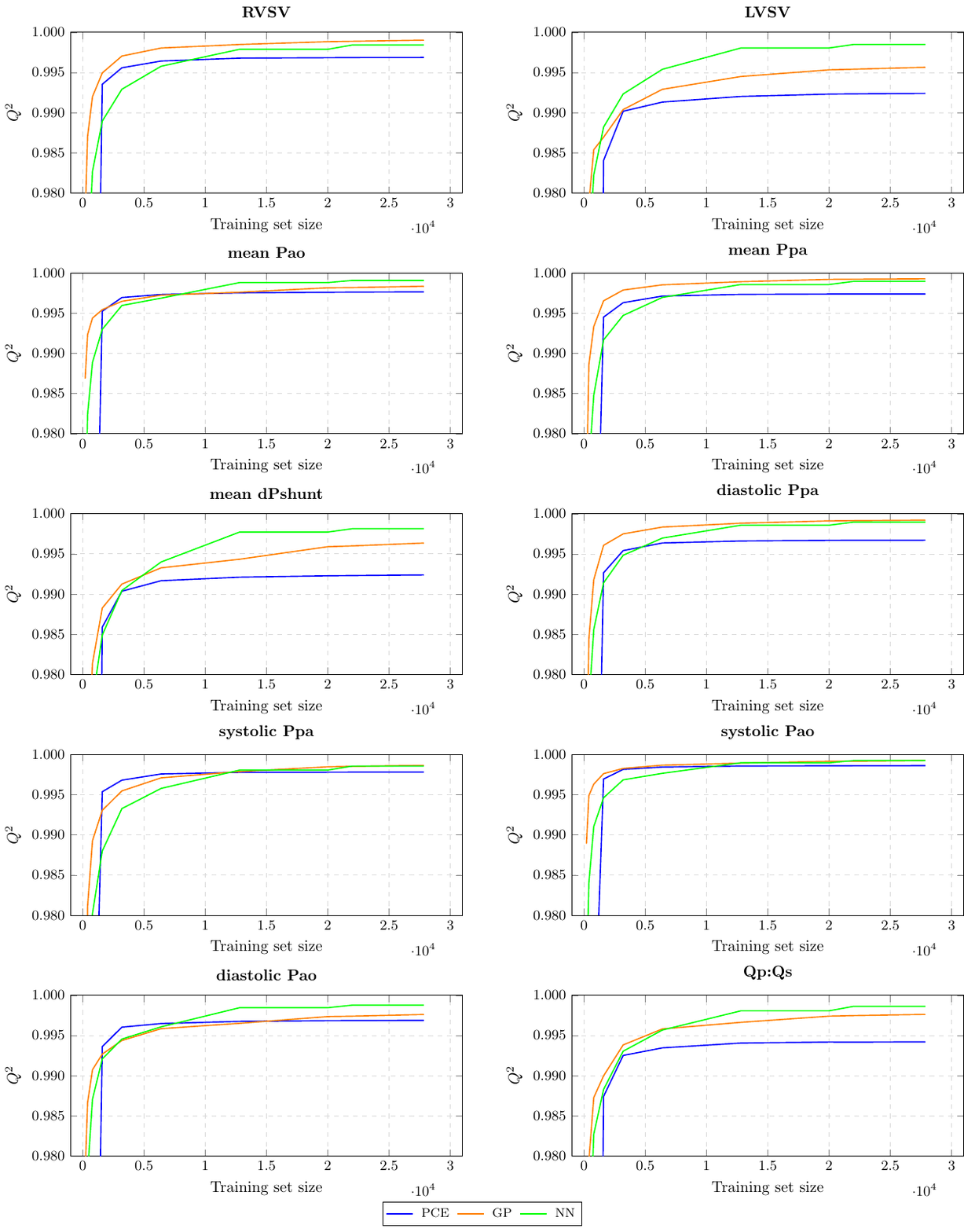} %.png or .pdf for latex
    \caption{$Q^2$ metric using the testing dataset for the second model.}
    \label{fig:q2_pavlos}
\end{figure}

The Q2 for all the outputs is plotted in figure~\ref{fig:q2_pavlos}.

\begin{figure}[H]
    \centering
    \includegraphics[width=0.9\columnwidth]{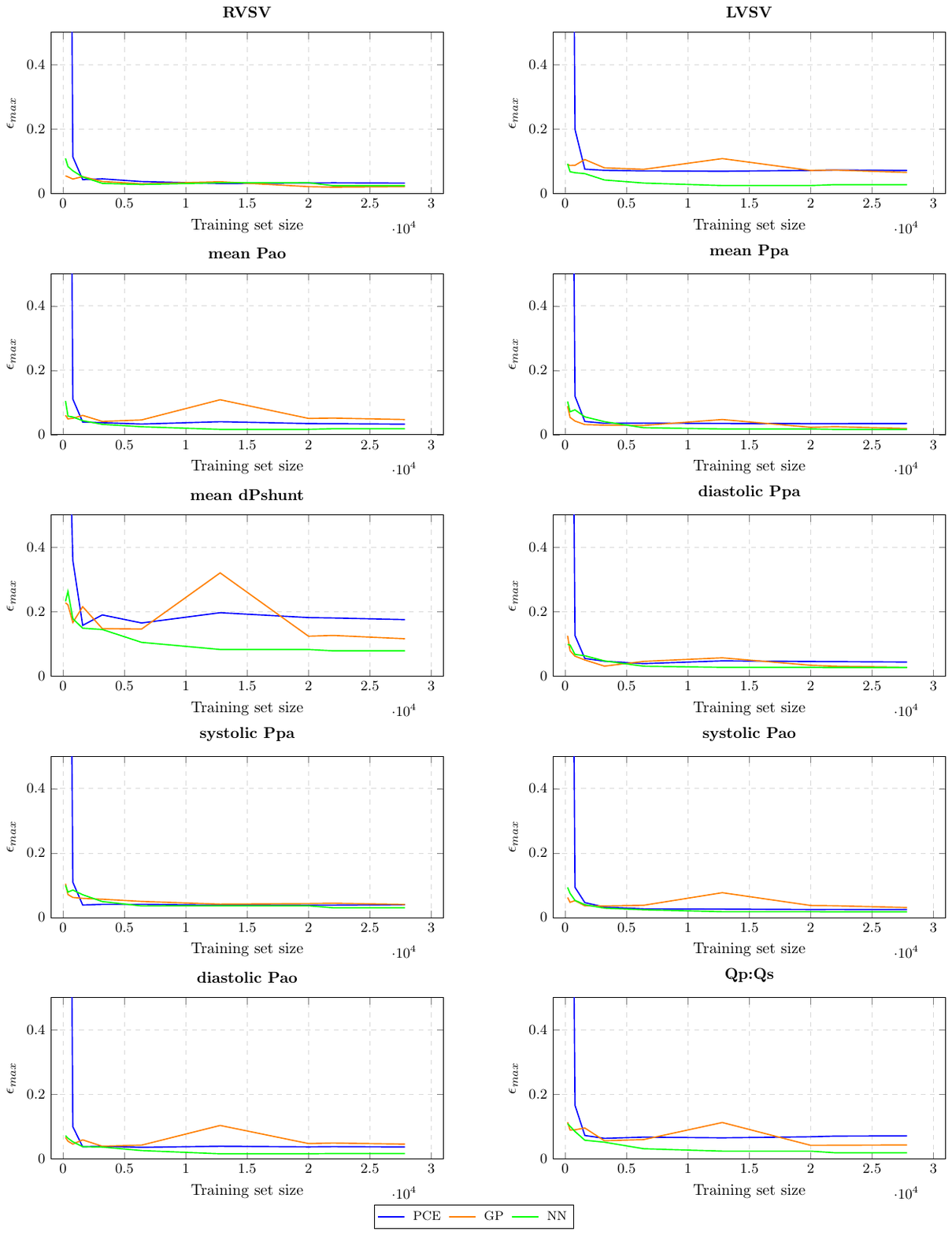} % .png or .pdf for latex
    \caption{Maximum relative error using the testing dataset for the second model.}
    \label{fig:max_err_pavlos}
\end{figure}

It can be noticed from Figure~\ref{fig:q2_pavlos} that neural networks outperform GP and PCE for large training datasets. While for quite small training sets GP performs the best, while PCE has poor performance. However, it should be noted that all methods have a Q2 value of above 0.98 for large training sets. To better compare the methods, the maximum absolute error is plotted against the training set size in figure~\ref{fig:max_err_pavlos}.

It can be seen from Figure~\ref{fig:max_err_pavlos} that neural networks have the lowest max relative error in the test set. The error of PCE and GP can go up to 10-20\% for certain outputs (Mean shunt pressure gradient), while neural networks have a maximum error of less than 10\%.

\subsection{Model 3}

For the third model, the inputs and outputs are shown in table~\ref{tab:model3inputsoutputs}.

\begin{table}[H]
    \centering
    \caption{Input and output parameters for Model 3}
    \begin{tabular}{|c|c|}
    \hline
    \textbf{Inputs} & \textbf{Outputs} \\
    \hline
    Vascular anomaly intensity ($s_0$) & Contrast agent concentration signal ($c(t)$) \\
    Cardiac output (CO) &\\
    Mean arterial pressure (MAP) &\\
    Liver extraction rate ($\alpha$) &\\
    Proportion of blood going to the organ of interest ($x_\mathrm{organ}$)&\\
    Body weight (m) &\\
    Specific blood volume ($BV_\mathrm{spec}$) &\\
    Heart rate ($H_r$) & \\
    Main recirculation time ($\tau$) &\\

    \hline
   
    \hline
    \end{tabular}
    \label{tab:model3inputsoutputs}
\end{table}

Since the output is a time-dependent variable, a different metamodeling technique needs to be used. Fourier transform was initially our first candidate, however, the signals are not periodic and most them are early cut and most of the frequencies knowledge is lost. 

GP and PCE metamodeling approaches are not well-suited for time-dependent models because they do not account for temporal correlations between outputs at different time points. A naive approach would treat each time step as an independent output, significantly increasing the dimensionality and computational cost, making these methods inefficient. 

Additionally, GPs and PCEs scale poorly with high-dimensional outputs, leading to impractical runtimes. Instead, methods specifically designed for sequential data, like recurrent neural networks (RNNs) or variants as Long Short-Term Memory (LSTM), are better suited for capturing the dependencies in time-dependent models.

Therefore, specifically for this model, we implemented an LSTM architecture in Tensorflow. The LSTM architecture consists in a few-to-sequence model outputting a sequence of size 200. In the architecture, only one LSTM block was implemented since it was enough to produce accurate predictions.

\begin{figure}[H]
\centering
    \includegraphics[width=0.9\columnwidth]{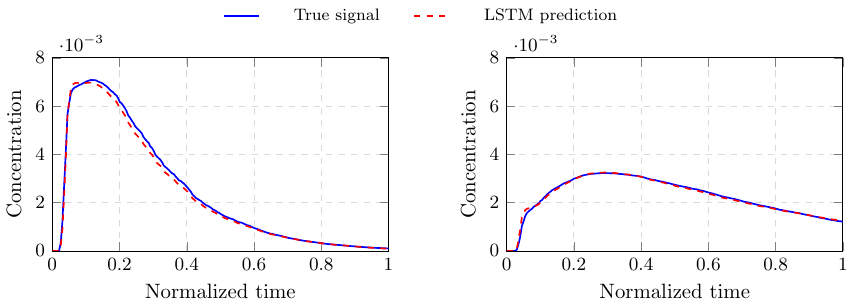}
\caption{Two testing samples of LSTM predictions of the contrast agent concentration. On the left, one can find a signal representing a healthy case (low vascular anomaly intensity $s_0$), while on the left a diseased case (higher vascular anomaly intensity $s_0$).}
\label{fig:rnn_results}
\end{figure}

Figure~\ref{fig:rnn_results} shows samples from the test sets and their corresponding LSTM predictions, which clearly shows good agreement. The two cases are chosen due to their different signal behavior corresponding to the diseased level of the patient.

\section{Applications of the metamodel}\label{section:application}

After the training phase of the metamodel, it can be used in several applications, including SA, parameter estimation, and UQ. In this section, we chose neural networks as a showcase since it was seen to outperform PCE and GP. The second 0D cardiovascular model is chosen for the application using neural networks.

\subsection{Sensitivity analysis}\label{section:SA}

Variance-based sensitivity analysis (Sobol indices), which is a global sensitivity analysis technique, is chosen in this work. When PCE is chosen as the surrogate model, Sobol indices can be derived directly from the coefficients of the polynomial terms \cite{sudret2008global}, whereas for NNs and GPs, a classical Monte Carlo approach is used \cite{sobol2001global}. In this section, we compare the total Sobol indices for the PAH model, first calculated for its 0D version using Saltelli's method, and then for its surrogate model built with neural networks.
%As will be discussed in Subsection 6.1, performing sensitivity analysis with meta-models is convenient for ensuring the physiological validity of the simulations. 

In order to conduct the sensitivity analysis for the 0D model, approximately 45,000 model evaluations were performed, reaching convergence. 

\begin{figure}[H]
\centering
\includegraphics[width=1.0\columnwidth]{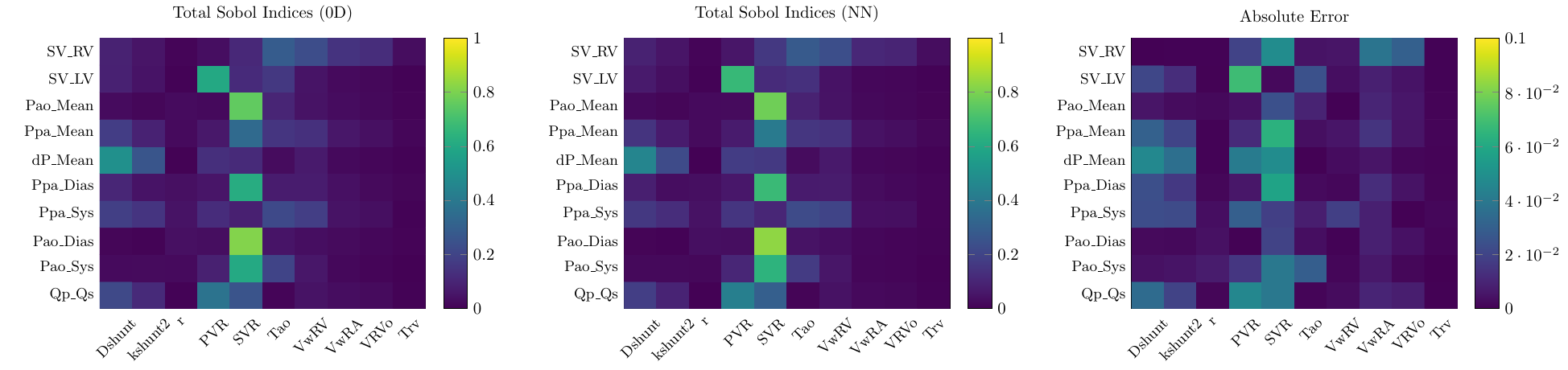}
\caption{Sobol indices for the second cardiovascular model. On the left, one can find the total Sobol indices using the 0D model, and on the middle using neural network predictions, and the absolute error is plotted on the right. }
\label{fig:test}
\end{figure}

According to figure \ref{fig:test}, both methods yield nearly identical results, with any differences being negligible. This sensitivity map indicates that patient-specific parameters, particularly Systemic and Pulmonary Vascular Resistances (SVR and PVR), as well as the contractility of the chambers ($T_{ao}$) and the wall volume of the RV ($V_{wRV}$), have the most significant impact on the various model outputs, and can therefore be considered as biomarkers for predicting surgery outcomes. %Additionally, it shows that shunt-related parameters, such as the diameter and shape of the shunt, also have a considerable impact on the outcomes and should be carefully considered in pre-operative planning.  

\subsection{Parameter estimation} \label{section:param_estimation}

Since the generated metamodel is differentiable and automatic differentiation is readily available in most deep learning frameworks, gradient descent-based optimization can be used to perform parameter estimation. Given measured output $\mathbf{\bar{Y}}$, the input parameters $\mathbf{X}$ can be estimated by solving an unconstrained minimization problem defined as:

\begin{equation}
\min_{\mathbf{X}} |\mathbf{\hat{Y}(\mathbf{X}; \hat{\theta}})-\mathbf{\bar{Y}}^2|,
\end{equation}

where $\hat{\theta}$ are the optimized neural network parameters (weights and biases), which are kept fixed.

Since the input parameters can have definitive ranges, a constrained form of the problem is defined as:

\begin{equation}
\min_{\boldsymbol{\psi}} |\mathbf{\hat{Y}(\boldsymbol{\sigma(\psi)}; \hat{\theta}})-\mathbf{\bar{Y}}^2|,
\end{equation}

where $\mathbf{X} = \boldsymbol{\sigma(\psi)}$, $\boldsymbol{\sigma}$ is a function designed to constrain to input space and is chosen as a variant of the sigmoid function $\mathbf{S}$, and $\boldsymbol{\psi}$ is a dummy variable. This construction is chosen to ensure that the estimated parameters falls within physiological ranges. The function $\boldsymbol{\sigma(\psi)}$ can be written as:

\begin{equation}
\boldsymbol{\sigma(\psi)}= \mathbf{X_1 + (X_2-X_1)S(\psi)},
\end{equation}

where $\mathbf{X_1}$ and $\mathbf{X_2}$ are the lower and upper bounds of the function.
It should be noted that this addition does not ensure the validity of the solution of the inverse problem since some of the input parameters might not be uniquely identifiable. However, if there is a single local minimum in that specified range, it helps reaching that specific minimum instead of falling in a non physiological minimum.

\begin{figure}[H]
    \centering
    \includegraphics[width=0.8\columnwidth]{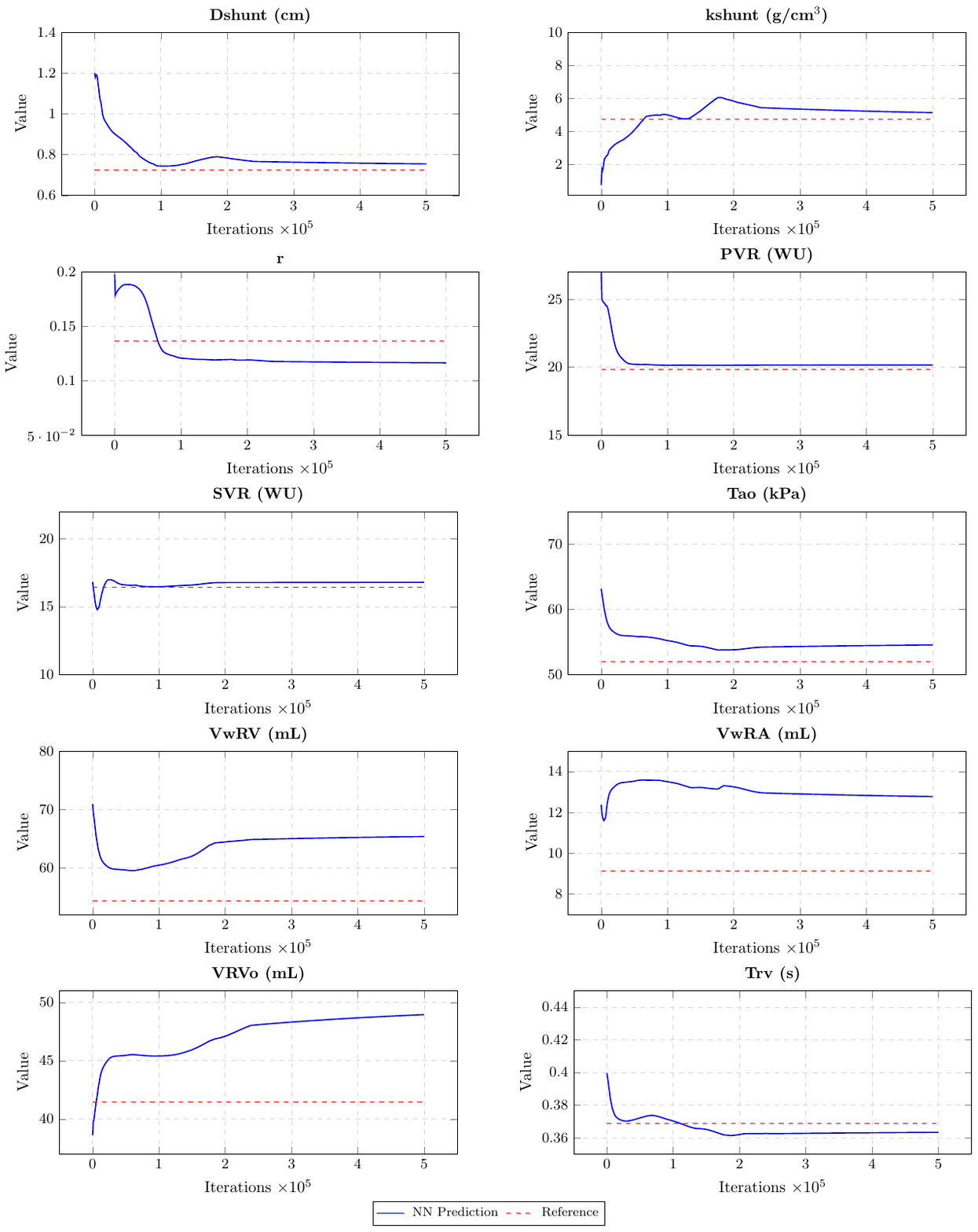} % .png or .pdf for latex
    \caption{The value of the input parameters vs. iterations showing the evolution of the parameters with the minimization algorithm steps.}
    \label{fig:inverse}
\end{figure}

To solve this minimization problem, $\boldsymbol{\psi}$ is initialized randomly and then updated through the minimization process with Adam optimization algorithm. An example of a solved inverse problem is shown in figure~\ref{fig:inverse} where 50000 Adam iterations are performed taking nearly 2 minutes. All parameter estimations converged.

It is noticed from figure~\ref{fig:inverse} that the parameters that are more sensitive, according to the Sobol indices results, were well identified as for PVR and SVR. However, parameters that are less sensitive were not identified accurately. That can mean that these parameters are not uniquely identifiable.

The error in the parameter estimation is shown in table~\ref{tab:error_inverse}, which is calculated as follows:

\begin{equation}
    Error(\%) = 100*\dfrac{X_{est} - X_{ref}}{X_{max}-X_{low}}
\end{equation}

where $X_{est}$ is the value of the estimated parameter, $X_{ref}$ the reference value, and $X_{max},\ X_{low}$ are the physiological limits of the parameter.

\begin{table}[htbp]
    \caption{Error in the parameter estimation problem with respect to the physiological ranges of each parameter.}
    \centering
    \begin{tabular}{|c|c|}
        \hline
        \textbf{Parameter} & \textbf{Error in \%} \\ \hline
        Dshunt & 3.6 \%  \\  \hline
        Kshunt & 3.8 \% \\ \hline
        r & 13.4 \% \\ \hline
        PVR & 2.7 \% \\ \hline 
        SVR & 2.8 \% \\ \hline
        Tao & 10.3 \% \\ \hline
        VwRV & 39.4 \% \\ \hline
        VwRA & 45.5 \% \\ \hline
        VRVo & 53.3 \% \\ \hline
        Trv & 5.4 \% \\ \hline
    \end{tabular}
    \label{tab:error_inverse}
\end{table}

\subsection{Uncertainty quantification}
\begin{figure}[H]
\centering
\begin{subfigure}{0.45\textwidth}
    \includegraphics[width=\textwidth]{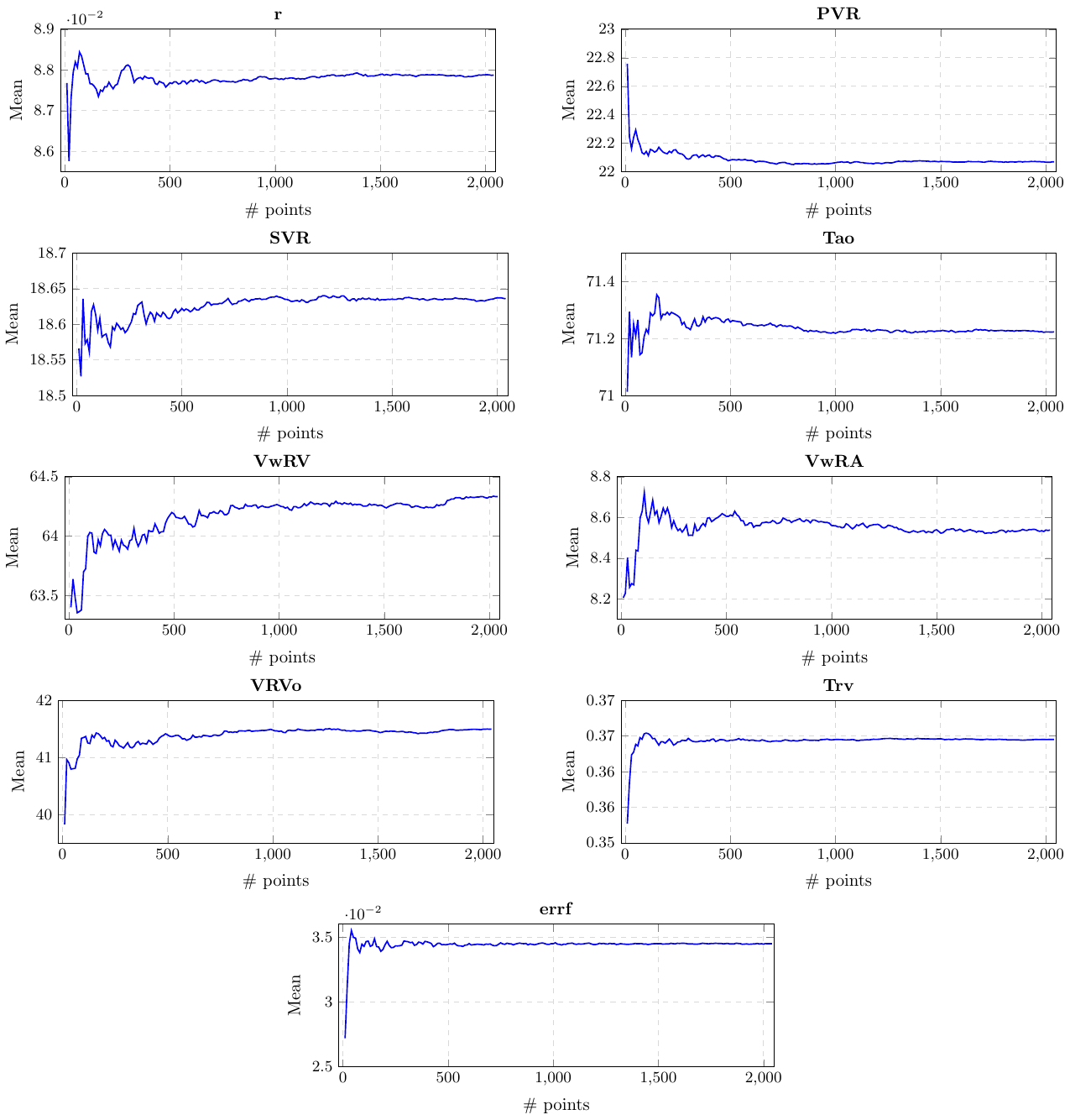}
    \caption{Convergence of input parameter mean as a function of the number of samples used}
    \label{fig:UQ_mean}
\end{subfigure}
\hfill
\begin{subfigure}{0.45\textwidth}
    \includegraphics[width=\textwidth]{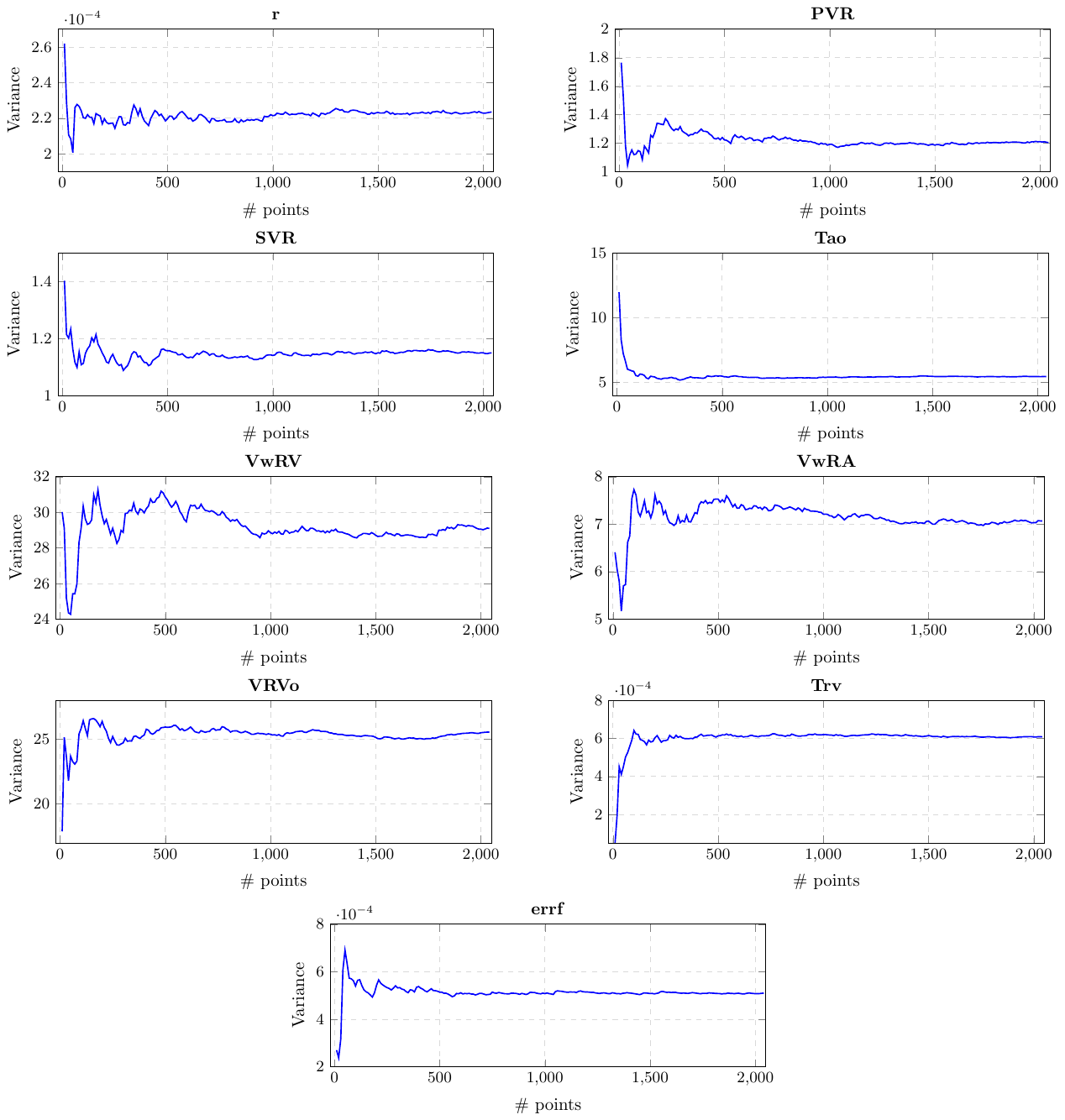}
    \caption{Convergence of input parameter variance as a function of the number of samples used}
    \label{fig:UQ_var}
\end{subfigure}
       
\caption{Convergence criteria of the inverse uncertainty quantification. As more samples are used for the inverse problem, mean and variance of each input parameter converge to a steady statistic.}
\label{fig:UQ_inputs}
\end{figure}

Due to the relatively low computational cost of the NN emulator to solve inverse problems and forward predictions, patient specific UQ was performed using a two-stepped, Monte Carlo type approach. First, an inverse UQ is performed, whereby patient data and their assumed uncertainty (Table \ref{tab:patdata}) is sampled, and the inverse problem (as described in \ref{section:param_estimation}) is solved for each individual sample. Note that the inverse problem, as visualized in Figure \ref{fig_pipeline}, is solved with a NN that is trained on the pre-intervention model, i.e. without the Potts shunt in place, since the model should be optimized to pre-intervention data to reflect the patient state before treatment. Convergence to a stable distribution of optimized parameters is monitored by plotting the mean and variance over the number of samples solved (Fig \ref{fig:UQ_inputs}). Once convergence is achieved, the resulting optimized parameter distributions are then augmented with the model parameters which describe the placement of the shunt, and are used to predict the hemodynamic outcome of the placement of a Potts shunt. A similar result could be achieved by sampling a wider input range for the parameters describing the shunt resistance, which would lead to neglible flow across the shunt. However, this would likely significantly increase the number of samples required to train the NN.

\begin{table}
    \caption{Measured data  for an example patient, as published in earlier work \cite{pant2022multiscale}. All measurement errors are assumed independent and uniformly distributed.}
    \centering
    \begin{tabular}{|c|c|c|c|}
     \hline
         \textbf{Data type} & \textbf{Value} & \textbf{Unit} & \textbf{Assumed Measurement error} \\ \hline
         $P_{ao,dias}$ & 53 & mmHg & 2.5\% \cite{edwardsmanual} \\  \hline
         $P_{ao,mean}$ & 69 & mmHg & 2.5 \% \cite{edwardsmanual}\\ \hline
         $P_{ao,sys}$ & 94 & mmHg & 2.5 \% \cite{edwardsmanual}\\ \hline
         $P_{pa,dias}$ & 67 & mmHg & 2.5 \% \cite{edwardsmanual}\\ \hline
         $P_{pa,mean}$ & 85 & mmHg & 2.5 \% \cite{edwardsmanual}\\ \hline
         $P_{pa,sys}$ & 112 & mmHg & 2.5 \% \cite{edwardsmanual}\\ \hline
         $SV_{LV}$ & 51 & mL & 8 \% \cite{mercado2017transthoracic} \\ \hline
         $SV_{RV}$ & 51 & mL & 8 \% \cite{mercado2017transthoracic} \\ \hline
    \end{tabular}

    \label{tab:patdata}
\end{table}

\begin{figure}[H]
    \centering
    \includegraphics[width=0.75\linewidth]{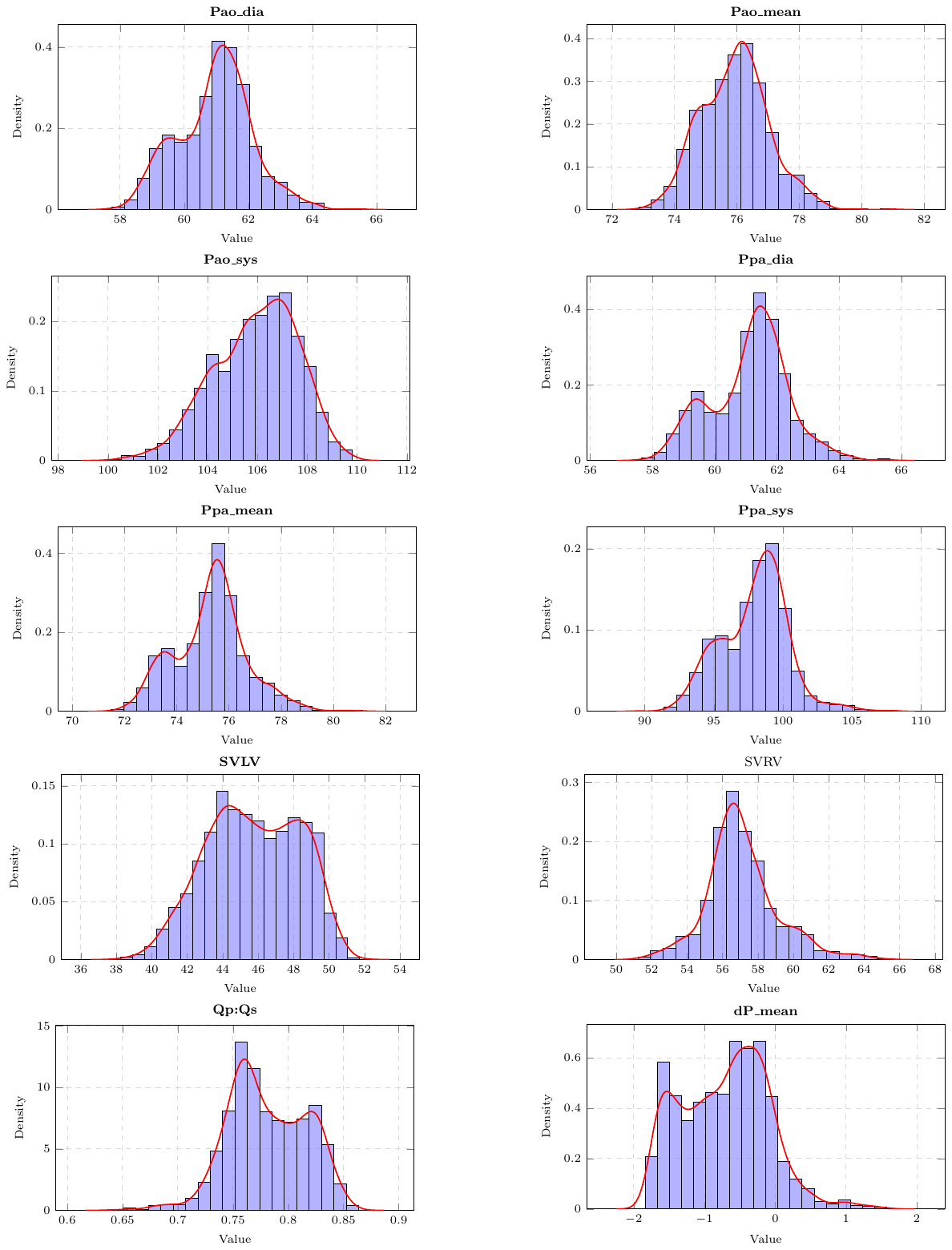}
    \caption{Histograms indicating the uncertainty in the patient-specific predicted values of a pulmonary hypertensive patient after placement of a Potts shunt, propagated from the uncertainty in the input parameters of the model before the shunt placement.}
    \label{fig:UQ_predictions}
\end{figure}

We demonstrate the UQ part of the pipeline on a pulmonary hypertensive patient presented in earlier work \cite{pant2022multiscale}, whose measured characteristics are presented in Table \ref{tab:patdata}. All measurement errors are assumed independent and uniformly distributed. Sampling was performed using a Sobol` sequence, and after 2048 samples parameter input means and variances were deemed converged \ref{fig:UQ_inputs}. The achieved samples are then used to solve the forward UQ problem, assuming a shunt placed with the following descriptive parameters of $D_{shunt}\  =\  1.0\ cm$ and $k_{shunt}\ =\ 3.0\ mL^2/(mmHg \cdot s^2$). Uncertain predictions are shown in the form of histograms in Figure \ref{fig:UQ_predictions}. Note that the predicted outcomes are no longer uniformly distributed, which is due to both the non-linearity of the model input-output relationships, as well as the additional step of solving the inverse problem. 

\section{Discussion and conclusion}\label{section:conclusion}

%From the obtained results of the different models, it can be deduced that PCE is outperformed when the training dataset is small (less than 800 points), while NN and GP provide acceptable results for such low dataset sizes.
In this section, several issues are discussed including computational time, ensuring the physiological validity of metamodel outputs, and performing inverse problems.

\subsection{Computational time}\label{sec:computational_time}

The computational time of the different metamodels employed is discussed in this section as shown in Table~\ref{t:cost}. For the building of metamodels, it was observed that neural networks have the fastest training time in the case of a big dataset ($> 30K$ data points). The computational time was always in the order of few minutes when training with the full dataset. However, this cost goes to nearly an hour for PCE and around 12 hours for GPs, when running on the same laptop. This is the price paid by GPs in order to offer a natural way to quantify uncertainty in predictions, providing the predicted mean but also a confidence interval, representing the uncertainty at each prediction point. This is particularly valuable in fields like Bayesian optimization, where understanding uncertainty is crucial. It should be noted that when training on a small part of the dataset, 3K points for example, GPs and PCE have training times in the order of minutes. 

Once an emulator is trained, obtaining Sobol indices using the Saltelli method \cite{saltelli2004sensitivity} is a relatively fast process for NN and GPs (few minutes), due to the low computational cost per emulator evaluation. PCE has the advantage that Sobol indices are readily available when the model has been trained, with no additional computational cost.

For the inverse problem solving and uncertainty quantification, we use TensorFlow's automatic differentiation to compute gradients for the optimizer in NN models. This process generally converges within a minute. While automatic differentiation could also be applied to GPs or PCEs, these models must first be implemented within TensorFlow or another automatic differentiation framework, which is not as straightforward for GPs and PCEs as it is for NNs.

A summary of the computational cost for all the metamodels is shown in table~\ref{t:cost}.

\begin{table}[H]
\caption{Computational cost of the different tasks for the three metamodels: neural networks, polynomial chaos expansion, and gaussian processes regression.  'Requires additional setup' indicates tasks that need extra implementation effort.}
\centering
\begin{tabular}{|c|c|c|c|}
\hline
\rowcolor{gray!20}  % Color the whole row
 & \textbf{Neural Networks} & \textbf{Polynomial Chaos expansion} & \textbf{Gaussian processes} \\
 
\textbf{Metamodel building} & \cellcolor{green!30} Low & \cellcolor{orange!30} Medium & \cellcolor{red!30} High\\

\textbf{Sensitivity Analysis} & \cellcolor{orange!30} Minutes & \cellcolor{green!30} Free & \cellcolor{orange!30} Minutes \\ 

\textbf{Inverse problem} & \cellcolor{green!30} $<$ A minute & \cellcolor{orange!30} Requires additional setup & \cellcolor{orange!30} Requires additional setup \\

\textbf{Uncertainty quantification } & \cellcolor{green!30} Low & \cellcolor{orange!30} Requires additional setup & \cellcolor{orange!30} Requires additional setup \\
\hline
\end{tabular}
\label{t:cost}
\end{table}

\subsection{Ensuring physiological validity in metamodels}
When developing metamodels for physiological systems, ensuring that the outputs adhere to physiological constraints is a fundamental challenge. This issue is prevalent in both full models and metamodels, especially when the training data includes instances that fall outside physiological ranges.
To address this, specific strategies, depending on the application, are needed to maintain the validity of results within physiological ranges.
The applications that we are interested in this paper are prediction, sensitivity analysis (SA), parameter estimation, and uncertainty quantification (UQ).

For prediction purposes, the challenge is to ensure that the outputs of the metamodel adhere to physiological constraints. This can be approached in two main ways. One method involves training the metamodel exclusively on data that is known to fall within physiological ranges. By doing so, the model learns to generate predictions that are inherently physiologically plausible. In our models, we chose to use known ranges of (patho-)physiological behaviour (i.e., pressure, flows and volumes in respective model components) to determine ‘physiological behaviour’, and filter out those samples which are outside of that behaviour \cite{sala2023sensitivity, Haghebaert2025}. However, this approach limits the model's applicability to the specific ranges of inputs and outputs included in the training data. The second method involves using post-prediction filtering. In this approach, the metamodel generates predictions based on a broader set of data, but these predictions are subsequently checked and corrected to ensure they fall within physiological limits. While this method helps ensure that final outputs are valid, it does not prevent the generation of non-physiological outputs during the modeling process itself.

In the context of SA, which aims to understand how variations in input parameters affect output variability, the approach can differ based on the tools used. 
One effective strategy proposed in \cite{sala2023sensitivity} is to build a metamodel based on PCE using a training set restricted to input-output pairs that are physiologically valid. 
This approach ensures that the PCE provides Sobol indices naturally, reflecting the sensitivity of outputs to inputs while avoiding biases from non-physiological data. 
For those who may not have access to PCE or are unfamiliar with it, another approach is to create a general metamodel, such as one based on GP or NN, using a dataset that includes both physiological and non-physiological pairs. Sobol indices can then be computed using the classical Monte-Carlo approach as shown in Sec. \ref{section:SA}.
%variance-based methods combined with a Markov Chain Monte Carlo (MCMC) approach. During this process, outputs generated through sampling that fall outside physiological ranges are discarded, ensuring that only valid results contribute to the sensitivity analysis and the computation of the variances for the SA \cite{sobol1993sensitivity}.

Parameter estimation presents another scenario where physiological constraints are important. When a metamodel is trained on a comprehensive dataset that includes both physiological and non-physiological outputs, the parameter estimation process can be conducted without bias. The metamodel will provide accurate predictions across the input parameter space, and any non-physiological outputs encountered during estimation can be identified and discarded. This method ensures that the parameter estimation algorithm operates without being constrained by a restricted training dataset. Conversely, if the metamodel has been trained solely on physiological data, parameter estimation can become more complex. The estimation algorithm might encounter difficulties if it needs to explore input parameter spaces that were not covered during training. This issue can be particularly problematic if the relationship between inputs and outputs is nonlinear and complex; for instance, two or more outputs are dependent on each other in generating non-physiological outputs. In such cases, the parameter estimation process must be adapted to account for the constrained training data and ensure that the search for valid input parameters remains effective.

Finally, UQ involves assessing how uncertainty in input parameters affects outputs. In principle, UQ faces similar challenges as SA with respect to maintaining physiological constraints. However, in practice, the narrower ranges of input variations typically encountered in UQ studies often do not result in issues with non-physiological outputs. This is because the range of input variations in assessing uncertainty for a given patient is usually smaller compared to the broader sampling over an entire population required for SA, thus reducing the likelihood of generating non-physiological results.

In summary, ensuring physiological validity in metamodel outputs requires tailored strategies depending on the application. For prediction, either training on physiological data or applying post-prediction filtering can be used. SA can benefit from PCE-based metamodels or general metamodels with filtering techniques. Parameter estimation needs careful consideration of whether the metamodel was trained on all data or only physiological pairs. Lastly, UQ generally handles narrower input ranges, mitigating issues with non-physiological outputs. 
By employing the proposed strategies appropriately, metamodels can provide reliable and physiologically relevant results across various applications.

\subsection{Parameter estimation and inverse problems}

The choice of the inputs to the neural network are typically the most impactful input parameters of the 0D models. This can be decided with clinicians in order to have the most clinically relevant parameters (that can possibly serve as biomarkers) as inputs, or through performing sensitivity analysis as performed in the second model (Sec. \ref{sec:model_2}). The outputs of the neural networks are typically clinically relevant outputs of the 0D model.  

Some of these outputs are clinically measured, while input parameters need to be estimated. For instance, the shunt diameter may be determined based on the clinically preferred stroke volume and ventricular pressures. This is achieved after building the surrogate model by solving an inverse problem, as demonstrated in section~\ref{section:application}.

While outside of the scope of this work, in some cases there are several values or combinations of inputs that lead to the same output, \textit{i.e.}, the chosen set of input parameters is unidentifiable. Depending on the research question, there are several approaches to deal with unidentifiability. First is a Monte Carlo type approach, where the inverse problem is solved several times with different random initializations, leading to a distribution of input parameters \cite{shin2019scalable}. Some methods leverage the sensitivity indices to derive an orthogonal set of input parameters \cite{li2004selection, saxton2023personalised}, while other methods allow for a parametric description of non-identifiable parameter combinations \cite{tong2024invaert}.
%, Tonkin}.
Alternatively, a formal identifiability analysis is performed, iteratively reducing the input space by fixing unidentifiable parameters to produce a reduced, identifiable set; such as profile likelihood analysis \cite{maiwald2016driving, vanlier2012integrated}, diaphony \cite{hornfeck2015diaphony, van2020parameter} or determining intra-class correlation \cite{koopsen2024parameter}. However, these methods can be prohibitive due to their computational cost, which can be overcome by using surrogate models.

\subsection{Conclusion}

In this study, various metamodeling strategies are evaluated on three different 0D models.

Our findings show that PCE does not perform well with small training datasets but excels with larger ones. A key advantage of PCE is its ability to derive Sobol indices without additional computational cost.

GPs offer a natural way to quantify uncertainty by providing both a predicted mean and confidence intervals, which is crucial, for instance, in Bayesian optimization. As non-parametric models, GPs adjust their complexity based on the data, offering flexibility, and they are data-efficient, performing well with small datasets due to their probabilistic nature. However, GPs do not scale well with large datasets and become computationally expensive, in contrast to Neural Networks (NNs), which are better suited for large-scale problems due to their scalability.

NNs, on the other hand, offer scalability to large datasets, with relatively low computational cost and high accuracy. They are versatile and straightforward for solving inverse problems and uncertainty quantification, particularly due to the automatic differentiation capabilities available in libraries such as TensorFlow. Additionally, architectures like Long Short-Term Memory (LSTM) networks handle time-series data effectively, providing strong approximation capabilities. In contrast, both PCE and GPs struggle to efficiently approximate time-series data, and moreover, their computational cost grows significantly, making them less suited for such tasks.

%We demonstrated how NN could effectively represent a variety of cardiovascular models, of different complexities, and types of outputs (varying-in-time or not).

To conclude, overall NN present advantages in terms of implementation and run-time cost-effectiveness to perform the whole pipeline of studying and performing patient-specific modeling.  

\section*{Code and data availability}
Codes and data are available on GitLab at \url{https://gitlab.com/0d-cardiovascular/metamodeling}.

\section*{Acknowledgments}
We acknowledge funding from the European Research Council (ERC) under the European Union’s Horizon 2020 research and innovation program (Grant agreement No. 864313). 

\bibliography{biblio}

\appendix
\section{Estimation of computational time for full 0D models}

In this appendix, we present the estimated computational time for the different tasks using the full 0D models as shown in the table below.

\begin{table}[H]
\centering
\caption{Computational cost of different tasks for three models using the original 0D model.}
\begin{tabularx}{\textwidth}{|X|X|X|X|X|X|}
\hline
\rowcolor{gray!20}
\centering \textbf{Model/Task} & \centering \textbf{1 Simulation} & \centering \textbf{Sensitivity Analysis} & \centering \textbf{Inverse Problem} & \centering  \textbf{Uncertainty Quantification (UQ)} & \textbf{UQ of Inverse Problem} \\
\hline
\hline
\centering \textbf{Model 1} &\centering  1.93s in C & \centering $>$3 days &\centering  $\sim$20 mins & \centering $\sim$8h & $\sim$97 days \\
\hline
\centering \textbf{Model 2} &\centering  0.3s in Fortran &\centering  $\sim$4h &\centering  $\sim$5 mins &\centering  $\sim$50 mins & $\sim$34 days \\
\hline
\centering \textbf{Model 3} &\centering  5.6s in Python &\centering  $\sim$2 days &\centering  $\sim$5h &\centering  $\sim$15.5h & $\sim$59 weeks \\
\hline
\end{tabularx}
\label{t:cost_models}
\end{table}

\end{document}